\documentclass[12pt]{article}
\usepackage{}
\usepackage{amssymb}
\usepackage{amsmath}
\usepackage{amsthm}
\usepackage{amsfonts}
\usepackage{mathrsfs}
\usepackage{color, dsfont}
\usepackage[T1]{fontenc}
\allowdisplaybreaks
\usepackage{bm}
\usepackage{verbatim}
\usepackage{txfonts, graphicx, nicefrac}

\usepackage[noend]{algpseudocode}
\usepackage{algorithmicx}
\usepackage{algorithm,float}

\allowdisplaybreaks

\def\az{\alpha}  \def\bz{\beta}
    \def\dz{\delta}
    
\def\gz{\gamma}  \def\kz{\kappa}
\def\lz{\lambda} \def\mz{\mu}
\def\nz{\nu}     \def\oz{\omega}
\def\pz{\pi}

\def\rz{\rho}        \def\sz{\sigma}
        
\def\vz{\varepsilon} \def\xz{\xi}

\def\llz{\Lambda} \def\ddz{\Delta}
    
\def\ggz{\Gamma}  \def\ooz{\Omega}
\def\ppz{\Pi}     \def\ssz{\Sigma}

   \def\bq{{\mathscr{B}}}
   
   \def\fq{{\mathscr{F}}}

   \def\llq{{\mathscr{L}}}
   
   \def\pq{{\mathscr{P}}}

\def\qd{\quad}
\def\qqd{\qquad}

\def\lt{\left}
\def\rt{\right}

\def\PP{\mathbb{P}}
\def\EE{\mathbb{E}}

\newcommand{\mathsym}[1]{{}}

\def\leq{\leqslant}
\def\geq{\geqslant}

\newtheorem{thm}{Theorem}[section]
\newtheorem{prop}[thm]{Proposition}
\newtheorem{lem}[thm]{Lemma}

\newtheorem{cor}[thm]{Corollary}
\newtheorem{defn}[thm]{Definition}

\newtheorem{ass}[thm]{Assumption}

\numberwithin{equation}{section} \allowdisplaybreaks[4]

\makeatletter
\newenvironment{breakablealgorithm}
  {
   \begin{flushleft}
   \allowdisplaybreaks
     \refstepcounter{algorithm}
     \hrule height.8pt depth0pt \kern2pt
     \renewcommand{\caption}[2][\relax]{
       {\raggedright\textbf{\ALG@name~\thealgorithm} ##2\par}%
       \ifx\relax##1\relax 
         \addcontentsline{loa}{algorithm}{\protect\numberline{\thealgorithm}##2}%
       \else 
         \addcontentsline{loa}{algorithm}{\protect\numberline{\thealgorithm}##1}%
       \fi
       \kern2pt\hrule\kern2pt
     }
  }{
     \kern2pt\hrule\relax
   \end{flushleft}
  }
\makeatother

\def\prf{\medskip \noindent {\bf Proof}. }
\def\deprf{\quad $\square$ \medskip}

\def\be{\begin{equation}}
\def\de{\end{equation}}

\def\dear{\end{eqnarray*}}
\def\lb{\label}

\def\dps{\displaystyle}

\def\den{\end{enumerate}}

\def\d{\mathrm{d}}

\def\Ho{\mathcal{H}}

\def\So{S}
\def\Ao{A}
\def\Bo{B}
\def\Ko{K}

\def\PP{\mathbb{P}}
\def\EE{\mathbb{E}}

\def\RR{\mathbb{R}}

\def\lan{\langle}
\def\ran{\rangle}

\def\pla{{\it Player }}

\def\one{\mathds{1}}

\begin{document}
		\date{}
	\pagestyle{plain}
	\title{Discounted semi-Markov games with incomplete information on one side\footnote{{\bf Funding:} This work was partly supported by the National Natural Science Foundation of China (No. 11931018, 61773411, 11961005, 11701588), the Guangdong Basic and Applied Basic Research Foundation (No. 2020B1515310021) and the Natural Science Foundation of Guangdong Province (No. 2021A1515010057)}}	
	\author{Fang Chen \footnote{School of Mathematics, Sun Yat-Sen University, Guangzhou 510275, China. Email: chenf76@mail2. sysu.edu.cn}, Xianping Guo \footnote{School of Mathematics, Sun Yat-Sen University, Guangzhou 510275, China. Email: mcsgxp@mail. sysu.edu.cn}, Zhong-Wei Liao \footnote{corresponding author. College of Education for the Future, Beijing Normal University, Zhuhai 519087, China. Email: zhwliao@hotmail.com}}
	\date{}
	\maketitle \underline{}
	
{\bf Abstract:} This work considers two-player zero-sum semi-Markov games with incomplete information on one side and perfect observation. At the beginning, the system selects a game type according to a given probability distribution and informs to \pla 1 only. After each stage, the actions chosen are observed by both players before proceeding to the next stage. Firstly, we show the existence of the value function under the expected discount criterion and the optimality equation. Secondly, the existence and iterative algorithm of the optimal policy for \pla 1 are introduced through the optimality equation of value function. Moreove, About the optimal policy for the uninformed \pla 2, we define the auxiliary dual games and construct a new optimality equation for the value function in the dual games, which implies the existence of the optimal policy for \pla 2 in the dual game. Finally, the existence and iterative algorithm of the optimal policy for \pla 2 in the original game is given by the results of the dual game.
	
\vskip 0.2 in \noindent{\bf Key Words.} semi-Markov games, incomplete information, optimality equation, optimal policy, iterative algorithm
\vskip 0.2 in \noindent {\bf MSC 2020 Subject Classification.} Primary: 91A27, 91A35; secondary 91A25, 90C40, 93E20
	
\setlength{\baselineskip}{0.25in}
	
\section{Introduction}
The research of incomplete information repeated games originated from Aumann and Maschler in the 1960s (see \cite{AM68, AM95} and references therein). Many practical interactions are characterized by information asymmetries, the players are not fully informed the games information, such as the details of the system and their opponents actions. Therefore, the study of the game with incomplete information naturally arises. What Aumann and Maschler concerned was two-player zero-sum repeated games with lack of information on one side and perfect observation, which means that some system information or payoff functions are announced to Player 1 only. This model leads to novel strategic issues that cannot be adequately analyzed by focusing separately on either the uncertainty or the long-run aspect.
	
The games with incomplete information have a wide application prospect, which leads to a lot of works dedicated to generalizations of this model or close extensions of it. Let us cite for the advances of games with incomplete information on one side. Harsanyi \cite{H67} introduced one-stage Bayesian games with incomplete information. Renault \cite{R06} showed the existence of the value function in Markov chain games with incomplete information on one side. The optimal policy typically involves a repeated revelation of information because the state of system changes over time, see also H\"{o}rner et al. \cite{HRSV10}. Neyman \cite{N08} gave proved the existence of the value function and the optimal policy for all players of the repeated games, where each state follows a Markov chain independently of actions and at the beginning of each stage only \pla 1 is informed about the information. To analyze the optimal policy for players with missing information, the duality method is proposed by De Meyer \& Rosenberg \cite{MR99} and Laraki \cite{L02}. Cardaliaguet \cite{CP07}, Cardaliaguet \& Rainer \cite{CR09a, CR09b} studied the games with incomplete information, where the state variable was assumed to evolve according to a stochastic differential equation, see also Gr\"{u}n \cite{G12} and Oliu-Barton \cite{O15}. For the case of non-zero-sum games, the system will be more complex, and we refer to Hart \cite{H85}, Sorin \cite{S83}, Simon, Spie\.{z} \& Toru\'{n}czyk \cite{SST95} and Renault \cite{R00, R01} for further reading. It is worth noting that there are also a lot of results on the games with incomplete information on both sides, see Aumann \& Maschler \cite{AM68} and Gensbittel \& Renault \cite{GR15}. For more recent advances on the topic of incomplete information games, refer to Aumann \& Heifetz \cite{AH02} and Mertens, Sorin \& Zamir \cite{MSZ15}.
	
It is well known that the transition time of discrete-time Markov games is constant, while that of continuous-time Markov game satisfies the exponential distribution. However, in practical application, the transition  time may not satisfy either of these two situations. The semi-Markov game is a generalization of a discrete-time Markov game, where the transition time depends not only the present state and the actions chosen but also on the next state. In this sense, more general control problems come under the purview of the theory of semi-Markov processes than Markov processes, refer to Ja\'{s}kiewicz \cite{J02, J09}, Luque-V\'{a}squez \cite{L02}, Mondal \cite{M17}, Puterman \cite{P94} and Vega-Amaya \cite{V03} for further reading.
	
However,  to the best of our knowledge, there is no relevant research that considers the semi-Markov game with incomplete information. In this paper, we focus on the two-player zero-sum semi-Markov games with incomplete information on one side and perfect observation. In our models, the transition mechanism is defined by the semi-Markov kernels. The incomplete information on one side means that the system selects a game type $k \in K$ with probability $p_k \in [0, 1]$ at the beginning, but only announces this selection to \pla 1. Namely, \pla 2 is missing this information. After each stage, the actions chosen are perfectly observed by both players. The main feature of that interaction was when selecting an action now, one considers not only information revealed by ones action and its affect on future behavior of the opponent, but also its affect on the system. In the analysis of those games with incomplete information, the main difficulty is how to characterize the players' optimal policies. Especially, how the selection of the system affect the polices of the uninformed players and whether the uninformed players can conjecture the selection of the system through the actions of the opponents.

The main contributions of the present paper are as follows. 1) In our incomplete information model, the sequence of states follows the transfer mechanism of semi-Markov chain, which depends on the actions of the players. Our model generalizes the model of repeated games (see Aumann and Maschler \cite{AM68, AM95}) with incomplete information, which corresponds to the case identity transition matrix. Our results also generalizes the model of zero-sum semi-Markov games. Since semi-Markov games with incomplete information have a strong application prospect, it is of great significance to study the value function and the constructions of the optimal policy. 2) For different players, we propose different optimality equations, respectively. In particular, it is difficult to analyze the policy for the uninformed players by lack of information. We develop the dual game method proposed by De Meyer \cite{M96} and De Meyer \& Rosenberg \cite{MR99} (see also Laraki \cite{L02}), and then construct a new optimality equation which is different from that of the original games. 3) For the players' optimal policies, the conclusions we get are not only the existence, but also the effective and feasible iterative algorithms (see Algorithms \ref{algo-p1} and \ref{algo-p2}). In each step of the iterative algorithms, the selection mechanism $p \in \pq(K)$ of the original game $G(p)$ (or auxiliary vector $z \in \RR^{|K|}$ of the dual game $G^{\#} (z)$) needs to be recalculated according to the historical data. Therefore, the optimal policy is historical dependent rather than Markov. We fully believe that our analysis piques some theorists curiosity and pave the way toward a more complex model.
	
As described above, the paper is organized as follows. In Section 2, we formally introduce the model of semi-Markov games with incomplete information and the notations used, including a basic assumption and some preliminaries. Section 3 is devoted to proving the existence of the value function (Theorem \ref{s3-thm3}), introducing the optimality equation (Theorem \ref{OEthm}), and giving the algorithm for the value function (Corollary \ref{s3-Vn-prop}). The explicit construction of the optimal policy for \pla 1 (Theorem \ref{s4-thm1}) is proposed in Section 4, and the iterative algorithm (Algorithm \ref{algo-p1}) of optimal policies is also given. In Section 5, we introduce the concepts of dual games, which are connected to the original games by the dual variational formulas (Theorem \ref{s5-thm4}). We propose a new optimality equation of the value function in the dual game (Lemma \ref{s5-lem7}) and further show the existence of the optimal policy for \pla 2 in the dual game (Theorem \ref{s5-thm8}). Finally, a feasible iterative algorithm (Algorithm \ref{algo-p2}) of the optimal policy for \pla 2 is given in the end of this section.

\section{The game model and preliminaries}
{\bf Terminologies and notations.} In this paper, we adhere wherever possible to the following notations. Given a finite set $E$, denote by $|E|$ the cardinality of $E$. Let $\bq(E)$ be the Borel $\sz$-algebra of $E$ equipped with discrete topology, and $\pq(E)$ be the set of probability distributions on $E$. The set $\pq(E)$ is viewed as a subset of $\RR^{|E|}$ and for $p, q \in \RR^{|E|}$, denote by $\| p- q \| := \sum_{i\in E}|p_i - q_i|$. Hence, $(\pq(K), \| \cdot \|)$ becomes a complete separable metric space, see \cite{B68}. For any $i\in E$, $\dz_i$ denotes the Dirac measure on $i$. Write $\RR_+ := [0,\infty)$. Finally, if $X$ and $Y$ are Borel space, denote by $\pq(X|Y)$ the family of transition probabilities (or stochastic kernels) on $X$ given $Y$.

A zero-sum semi-Markov game with incomplete information model is defined by the collection:
\be\lb{game}
\{\Ko, \So, (\Ao \times \Bo), p, Q(\cdot, \cdot|i, a, b), c(k, i, a, b) \},
\de
where $\Ko$ is the set of game types, $\So$ is the set of states, and $\Ao$ and $\Bo$ are the set of actions of \pla 1 and \pla 2, respectively. $\Ko$, $\So$, $\Ao$ and $\Bo$ are assumed to be finite and equipped with Borel $\sz$-algebra. The probability $p \in \pq(K)$ is the law of the game type. Moreover, the transition mechanism of the semi-Markov game is defined by the semi-Markov kernel $Q(\cdot,\cdot|i,a,b)$ on $\RR_+ \times \So$ given $\So \times \Ao \times \Bo$, which is assumed that:
\begin{itemize}
	\item [\rm (i).] Given any $i, j\in \So$, $a \in \Ao$ and $b \in \Bo$, $Q(\cdot, j| i, a, b)$ is a nondecreasing and right continuous real-valued function on $\RR_+$ satisfying $Q(0, j|i, a, b)=0$;
	\item [\rm (ii).] For each $t \in \RR_+$, $Q(t,\cdot|\cdot,\cdot,\cdot)$ is a sub-stochastic kernel on $\So$ given $\So \times \Ao \times \Bo$;
	\item [\rm (iii).] The limit $ \lim_{t \rightarrow \infty} Q(t, \cdot | \cdot, \cdot, \cdot)$ is a stochastic kernel on $\So$ given $\So \times \Ao \times \Bo$, which means $\sum_{j \in \So}  \lim_{t \rightarrow \infty} Q(t,j |i, a, b) =1$ for all $(i, a, b) \in \So \times \Ao \times \Bo$.
\end{itemize}
If actions $a \in \Ao$ and $b \in \Bo$ are chosen at state $i$, then $Q(t,j | i,a,b)$ is the probability that the sojourn time in state $i$ is not greater than $t \in \RR_+$ and the system jumps into the next state $j$. It is possible that $j=i$ with positive probability. Finally, $c$ is a real valued function on $\Ko \times \So \times \Ao \times \Bo$ that denotes the payoff rate function and it represents the reward rate for \pla 1 and the cost rate for \pla 2. Since $\Ko$, $\So$, $\Ao$ and $\Bo$ are finite sets, the function $c$ is bounded. Without loss of generality, $c (k,i,a,b)$ is assumed to be nonnegative (equivalently, bounded below). For convenience, denote
$$
c^* := \max_{(k, i, a, b) \in \Ko \times \So \times \Ao \times \Bo} c (k, i, a, b).
$$

The rules of the semi-Markov game with incomplete information is as follows. At the beginning, one of game type $k$ is chosen according to $p \in \pq(K)$ and informed to \pla 1 only. At the initial decision epoch $t_0 = 0$, the system stays at $i_0 \in \So$, and \pla 1 chooses an action $a_0 \in \Ao$ according to the information of $i_0$ and $k$, meanwhile \pla 2 chooses an action $b_0 \in \Bo$ just according $i_0$. The actions $(a_0,b_0)$ are observed by both players. Then the system remains in $i_0$ for time $t_1$ and changes to $i_1$. At time $t_1$, \pla 1 selects an action $a_1 \in \Ao$ based on the game type $k$, $(i_0, a_0, b_0)$ and the current state $i_1$, while \pla 2 selects an action $b_1 \in \Bo$ based on $(i_0, a_0, b_0)$ and $i_1 $. Again, the actions $(a_1,b_1)$ are observed by both players. As a consequence of those action choices, the system remains in $i_1$ for time $t_2$ and changes to $i_2$. The game evolves repeatedly in the above way. The probability $p \in \pq(K)$ is a parameter in the analysis below, so the semi-Markov game with incomplete information is denoted by $G(p)$.

To introduce the history of the evolution of $G(p)$, we give the probability space (or trajectory space) which is based on Kitaev's construction \cite{K86, KR95}. Let $(\ooz, \fq)$ be the canonical measurable space that consists of the sample space
$$
\ooz := \Ko \times ( \So \times \Ao \times \Bo  \times \RR_+  )^{\infty},
$$
and the corresponding product $\sz$-algebra $\fq$. The elements of $\ooz$ are known as trajectories of the system. For each $\omega= (k,  i_0, a_0, b_0, t_1, \ldots, i_n, a_n, b_n,  t_{n+1}, \ldots) \in \ooz$, we define random variables $\kz$, $X_n$, $A_n$, $B_n$ and $T_n$ $(n=0,1,\ldots)$ on $(\ooz, \fq)$ as
\begin{align*}
& \kz (\omega) := k, \quad T_0(\omega):=0,\quad T_{n+1}(\omega) := \sum_{m=1}^{n+1} t_m, \\
& X_n (\omega) := i_n, \quad A_n (\omega) := a_n, \quad B_n (\omega) := b_n.
\end{align*}
The selection of game type is denoted by $\kz$ and the $n$-th decision epoch is $T_n$. $X_n$ is the state variable and $A_n$ (resp. $B_n$) is the action variable of \pla 1 (resp. \pla 2) at the $n$-th decision epoch. Hence, we define the processes of state and action until the $n$-th decision epoch by
\be\lb{h_n}
H_n(\omega):= (X_0,A_0,B_0,X_1,\ldots,A_{n-1},B_{n-1},X_{n})(\oz)=(i_0, a_0, b_0,i_1 \ldots,a_{n-1}, b_{n-1}, i_n).
\de
Let $\Ho_n = (\So \times \Ao \times \Bo)^{n} \times \So$ be the set of histories up to $n$-th decision epoch, which is equipped with its corresponding product $\sz$-algebra $\bq (\Ho_n)$. The finiteness of $\So$, $\Ao$ and $\Bo$ implies that $\Ho_n$ is finite too.

Next, we give the definitions of the policies for players. The asymmetry of informations between \pla 1 and \pla 2 leads to that the policies for \pla 1 depend on the game type $k \in \Ko$, but the policies for \pla 2 do not.

\begin{defn}\lb{pol}
\begin{itemize}
\item[(i).] A randomized history-dependent policy for \pla 1 is a sequence of stochastic kernels $\pz = \{ \pz^{(k)}_n, k \in \Ko, n \geq 0 \}$,  where $\pz^{(k)}_n$ is a stochastic kernel on $\Ao$ given $\Ho_n$, i.e.
		$$
		\pz^{(k)}_n (\cdot | h_n) \in \pq (\Ao), \qqd \forall h_n \in \Ho_n.
		$$
		Denote by $\sz = \{ \sz_n, n \geq 0 \}$, which is independent of game type $k$, the randomized history-dependent policy for \pla 2, where $\pz^{(k)}_n$ and $\Ao$ are replaced by $\sz_n$ and $\Bo$. Denote by $\ppz$ (resp. $\ssz$) the set of all randomized history-dependent policies for \pla 1 (resp. \pla 2).
		
		\item[(ii).] A deterministic policy for \pla 1 is a measurable mapping sequence $\varphi= \{ \varphi^{(k)}_n, k \in \Ko, n \geq 0 \}$, where $\varphi^{(k)}_n: \Ho_n \to \Ao$. Similarly, we can define the deterministic policy for \pla 2. Denote by $\ppz^D$ (resp. $\ssz^D$) the set of all deterministic policies for \pla 1 (resp. \pla 2).
	\end{itemize}
\end{defn}

Given any $p \in \pq(\Ko)$, $i \in \So$ and $(\pz, \sz) \in \ppz \times \ssz$, using Tulcea's Theorem (see \cite[Proposition C.10]{HL96}), there exists a unique probability measure $\PP_{p, i}^{\pz, \sz}$ on $(\ooz, \fq)$ such that
\begin{align}
\PP_{p, i}^{\pz,\sz} & (X_0=i, \kz=k) = p_k, \lb{sect2-0}\\
\PP_{p, i}^{\pz,\sz} & (A_n=a, B_n=b | \kz, H_{n}, T_n) = \pz^{(\kz)}_n ( a | H_n) \sz_n ( b | H_n), \lb{sect2-1} \\
\PP_{p, i}^{\pz,\sz} & (T_{n+1} - T_n \leq t, X_{{n+1}} = j | \kz, H_{n}, T_n, A_n, B_n) = Q(t, j| X_n, A_n, B_n). \lb{sect2-2}
\end{align}
Here and in what follow, we denote by $\EE_{p, i}^{\pz,\sz}$ the expectation with respect to $\PP_{p, i}^{\pz,\sz}$. According to (\ref{sect2-0})-(\ref{sect2-2}), for each measurable function $f$ on $\Ko \times (\RR_+ \times \So \times \Ao  \times \Bo)^{n+1} \times \RR_+ \times \So$, we obtain the expression of $\EE_{p, i}^{\pz,\sz} [f]$ as
\begin{align}\lb{E[f]}
& \EE_{p,i}^{\pz,\sz} \lt[ f( \kz, T_0, X_0, A_0, B_0, T_1-T_0,X_1, \ldots, A_n, B_n, T_{n+1}-T_n, X_{n+1}) \rt] \notag \\
& = \sum_{k \in \Ko} p_k \sum_{i_0 \in \So} \one_{\{i\}} (i_0) \sum_{a_0 \in \Ao, b_0 \in \Bo}\pz_0^{(k)} (a_0 | h_0) \sz_0 (b_0 |h_0 ) \sum_{i_1\in \So}\int_0^\infty Q(\d t_1, i_1 |i_0, a_0, b_0) \times \cdots \notag \\
& \qd \times \sum_{a_{n} \in \Ao, b_{n} \in \Bo} \pz_{n}^{(k)} (a_{n} |h_{n}) \sz_{n} (b_{n} | h_{n}) \sum_{i_{n+1}\in \So}\int_0^\infty Q (\d t_{n+1}, i_{n+1} | i_{n}, a_{n}, b_{n})\notag\\
&\qd \times f ( k, 0,i_0,a_0,b_0,t_1,i_1\ldots, a_n, b_n, t_{n+1}, i_{n+1}).
\end{align}
For each $p \in \pq (\Ko)$, we decompose that $p = \sum_{k \in \Ko} p_k \dz_{k}$, and then (\ref{E[f]}) implies that $\EE_{p, i}^{\pz, \sz}$ has the linearity property with respect to $p$, i.e.,
\be\lb{linear}
\EE_{p, i}^{\pz, \sz} [f] = \sum_{k \in \Ko} p_k \EE_{\dz_k, i}^{\pz, \sz} [f], \qd \forall i \in \So, \pz \in \ppz, \sz \in \ssz.
\de

Before introducing the goal of this paper, we give an assumption imposed on the semi-Markov kernel, in order to avoid infinite multiple decisions in a limited time. The following assumption is also used in \cite{HG10, K86, P94, R70}.
\begin{ass}\lb{ass1}
	There exist constants $\dz > 0$ and $0< \vz < 1 $ such that
	$$
	D (\dz | i, a, b) := \sum_{j \in \So} Q(\dz, j | i, a, b) \leq 1- \vz, \qqd \forall (i, a, b) \in \So \times \Ao \times \Bo.
	$$
\end{ass}

\begin{prop}\lb{s2-prop1}
Under Assumption \ref{ass1}, for each $p \in \pq(\Ko)$, $i \in \So$ and $(\pz, \sz) \in \ppz \times \ssz$, it holds that
$$
\PP_{p, i}^{\pz, \sz} \lt( T_{\infty} = \infty \rt) =1, \qd \text{where $T_{\infty}:= \lim_{n \rightarrow \infty} T_n$}.
$$
\end{prop}

\prf In fact, given any $M > 0$ and $n\geq 0$, we have
\be\lb{s2-prop2-1}
\PP_{p, i}^{\pz, \sz} \lt( T_{n+1} \leq M \rt) = \EE_{p, i}^{\pz, \sz} \lt[ \one_{T_{n+1} \leq M} \rt] \leq \EE_{p, i}^{\pz, \sz} \lt[ e^{- T_{n+1} + M} \rt], \qd \forall p \in \pq(\Ko), i \in \So, (\pz, \sz) \in \ppz \times \ssz.
\de
To calculate $\EE_{p, i}^{\pz, \sz} \lt[ e^{- T_{n+1}} \rt]$, we start with
$$
\EE_{p, i}^{\pz, \sz} \lt[ e^{- T_{n+1}} \rt] = \EE_{p, i}^{\pz, \sz}\lt[ e^{- T_{n}} \EE_{p, i}^{\pz, \sz} \lt(e^{- (T_{n+1} - T_{n})} \Big| H_{n}, T_{n},A_n,B_n \rt)\rt].
$$
According to (\ref{sect2-2}) and Assumption \ref{ass1}, we have
\begin{align*}
\EE_{p, i}^{\pz, \sz} \lt(e^{- (T_{n+1} - T_{n})} \Big| H_{n} , T_{n} ,A_n,B_n \rt) &= \int_0^{\infty} e^{-t} Q (\d t, S |X_n, A_n, B_n) \\
&= \int_0^\dz e^{-t} D(\d t | X_{n}, A_{n}, B_{n}) + \int_\dz^\infty e^{-t} D(\d t | X_{n}, A_{n}, B_{n}) \d t \\
& \leq (1- e^{-\dz}) D(\dz | X_{n}, A_{n}, B_{n}) + e^{- \dz} \leq 1 - \vz + \vz e^{- \dz} < 1.
\end{align*}
Using mathematical induction, we have
\be\lb{s2-prop2-2}
\EE_{p, i}^{\pz, \sz} \lt[ e^{- T_n} \rt] \leq (1 - \vz + \vz e^{- \dz})^n, \qd n \geq 1.
\de
Hence, by (\ref{s2-prop2-1}) we have $\PP_{p, i}^{\pz, \sz} \lt( T_{n+1} \leq M \rt) \leq e^M (1 - \vz + \vz e^{- \dz})^{n+1}$. Finally, the continuity of probability measure imply that $\PP_{p, i}^{\pz, \sz} ( T_\infty \leq M ) = \lim_{n \rightarrow \infty} \PP_{p, i}^{\pz, \sz} \lt( T_n \leq M \rt) = 0$. Combining with the arbitrariness of $M>0$, we have $\PP_{p, i}^{\pz, \sz} (T_\infty = \infty ) =1$.
\deprf

In order to characterize the value function of $G(p)$, we need to consider the continuous form of the stochastic processes $X_n$, $A_n$ and $B_n$. For each $\oz \in \ooz$ and $n \geq 0$, we define the processes at the time interval $t\in [T_n (\oz), T_{n+1} (\oz))$ as
\be\lb{pro}
X_t (\oz) = X_n (\oz) = i_n, \qd A_t (\oz) = A_n(\oz) = a_n, \qd B_t (\oz) = B_n (\oz) = b_n.
\de
Fix any discount factor $\az > 0$. Under Assumption \ref{ass1}, for any initial state $i \in \So$ and policies $\pz \in \ppz$ for \pla 1, $\sz \in \ssz$ for \pla 2, the expected discount reward for the \pla 1 of the game $G(p)$ is defined as
$$
V (p, i, \pz, \sz) := \EE_{p, i}^{\pz, \sz} \lt[ \int_0^{T_{\infty}} e^{-\az t} c (\kz, X_t, A_t, B_t) \d t \rt] = \EE_{p, i}^{\pz, \sz} \lt[ \int_0^\infty e^{-\az t} c (\kz, X_t, A_t, B_t) \d t \rt].
$$
To define our optimality criteria, we need to introduce the following concepts. The lower value of $G(p)$ is given by
\be\lb{lower}
\underline{V} (p, i) := \sup_{\pz \in \ppz} \inf_{\sz \in \ssz} V(p, i, \pz, \sz),
\de
which is called the game floor of \pla 1. Similarly, the upper value of $G(p)$ is given by
\be\lb{upper}
\overline{V} (p, i) := \inf_{\sz \in \ssz} \sup_{\pz \in \ppz} V(p, i, \pz, \sz),
\de
which is called the loss ceiling of \pla 2. It is clear that $\underline{V} \leq \overline{V}$. Conversely, if $\underline{V} \geq \overline{V} $, then we denote by $V^*$ the common value, which is called the value function of $G(p)$. The definition of the optimal policies is given below.

\begin{defn}\lb{optimal}
A policy $\pz^* \in \ppz$  for \pla 1 is called optimal in  $G(p)$ if
 $$\inf_{\sz \in \ssz} V(p, i, \pz^*, \sz) \geq \underline{V} (p, i), \qd \forall i \in S.$$
Similarly, a policy $\sz^* \in \ssz$  for \pla 2 is called optimal in $G(p)$ if $$\sup_{\pz \in \ppz} V(p, i, \pz, \sz^*) \leq \overline{V}  (p, i) \qd i \in S.$$
\end{defn}

The goal of this paper is to show the existence of the value function and find the optimal policies for \pla 1 and \pla 2 in $G(p)$. At the end of this section, we give the following result to explain the relationship of $\PP_{p, i}^{\pz, \sz}$ between different variables of $p \in \pq (\Ko)$,  $\pz\in\Pi$ and $\sz\in \ssz$.

\begin{prop}\lb{s3-thm5}
	Fix any $\lz \in [0, 1]$.
\begin{itemize}
\item [\rm i)]For each $p, q \in \pq (\Ko)$ and $\pz, \hat{\pz} \in \ppz$, there is $\pz^{\lz} \in \ppz$ such that
	\be\lb{s3-thm5-0}
	\lt( \lz \PP_{p, i}^{\pz, \sz} + (1- \lz) \PP_{q, i}^{\hat{\pz}, \sz} \rt) (E) = \PP_{\lz p + (1- \lz) q, i}^{\pz^{\lz}, \sz} (E) \qd \forall \, E \in \fq, \ i \in \So, \ \sz \in \ssz.
	\de
	Particularly, when $\pz=\hat{\pz}$, the above formula holds with $\pz^{\lz}= \pz$.
\item [\rm ii)] For each $\sz, \hat{\sz}\in \ssz$, there is $\sz^\lz$ such that
\be\lb{s3-thm5-0-0}
\lt( \lz \PP_{p, i}^{\pz, \sz} + (1- \lz) \PP_{p, i}^{\pz, \sz} \rt) (E) = \PP_{p, i}^{\pz, \sz^{\lz}} (E) \qd \forall \, E \in \fq,\ p\in \pq(K), \ i \in \So, \ \pz \in \ssz.
\de
\end{itemize}
\end{prop}
\prf \underline{Step 1, proof of (\ref{s3-thm5-0}).} Let $\{h_n \in \Ho_{n}, n\geq 0\}$ be the history sequence satisfying $h_0=i_0$ and $h_{n}=(h_{n-1}, a_{n-1}, b_{n-1}, i_{n})$. Given any $k \in \Ko$, define a sequence $\{C_n (k, h_n, a), n\geq 0,a \in A\}$ and $C_{-1}(k)$ as follow:
$$
\left\{
  \begin{array}{ll}
    C_{-1}(k):= \lz p_k + (1-\lz) q_k, &  \\
    C_n (k, h_n, a) := \lz p_k \pz_n^{(k)}(a|h_n) \dps\prod_{m=0}^{n-1} \pz_m^{(k)} (a_m | h_m)  + (1- \lz) q_k\hat{\pz}_n^{(k)}(a|h_n) \dps\prod_{M=0}^{n-1} \hat{\pz}_M^{(k)} (a_M | h_M). &
  \end{array}
\right.
$$
For convenience, use $C_{-1}(k,h_{-1},a_{-1})$ instead of $C_{-1}(k)$. $h_n=(h_{n-1},a_{n-1},b_{n-1},i_n)$ implies that $C_n (k, h_n, a) = 0$ for all $a \in A$ when $C_{n-1} (k, h_{n-1}, a_{n-1}) = 0$ and
\be\lb{s3-thm5-1}
\sum_{a \in \Ao} C_n (k, h_n, a) = C_{n-1} (k, h_{n-1}, a_{n-1}), \qd n \geq 0.
\de
Using $C_n (k, h_n, a)$, we define the policy $\pz^\lz = \lt\{ \pz_n^{\lz, (k)}, n \geq 0, k \in \Ko \rt\}$ as
\be\lb{s3-thm5-2}
\pz_n^{\lz, (k)}  (a | h_n):=
\left\{
\begin{array}{lll}
	C_n (k, h_n, a) / C_{n-1} (k, h_{n-1}, a_{n-1}), & \hbox{$C_{n-1} (k, h_{n-1}, a_{n-1}) > 0$;} \\
	| A |^{-1}, & \hbox{otherwise.}
\end{array}
\right.
\de
According to (\ref{s3-thm5-1}), for all $n \geq 0$ and $k \in \Ko$, it holds that $\sum_{a \in \Ao} \pz_n^{\lz, (k)}  (a | h_n)= 1$. Denote
\begin{align*}
E_{-1} &:= \{k\} \times ( \So \times \Ao \times \Bo\times    \RR_+ )^{\infty}, \\
E_n &:=  \{k\} \times \{(i_0, a_0, b_0)\}\times [0, s_1] \times \cdots \times \{(i_n, a_n, b_n)\}\times [0, s_{n+1}] \times( \So \times \Ao \times \Bo  \times \RR_+ )^{\infty},
\end{align*}
which are the measurable cylinder subsets of $(\ooz, \fq)$. Obviously, we have
$$
\lt( \lz \PP_{p, i}^{\pz, \sz} + (1- \lz) \PP_{q, i}^{\hat{\pz}, \sz} \rt) \lt( E_{-1} \rt) = \lz p_k + (1- \lz) q_k = \PP_{\lz p + (1- \lz) q}^{{\pz^\lz}, \sz} \lt( E_{-1} \rt).
$$
For $n \geq 0$, using (\ref{sect2-1}) and (\ref{sect2-2}), then
\begin{align}\lb{s3-thm5-3}
&\lt( \lz \PP_{p, i}^{\pz, \sz} + (1- \lz) \PP_{q, i}^{\hat{\pz}, \sz} \rt) \lt( E_n \rt) \notag \\
& = \Bigg\{ \lz p_k \one_{\{i\}}(i_0) \bigg[\prod_{m=0}^{n-1}\pz^{(k)}_m(a_m|h_m)\sz_{m}(b_m|h_m)\int_0^{s_{m+1}}Q( \d t, i_{m+1}|i_m,a_m,b_m)\bigg] \notag \\
& \qd \qd \times \bigg( \pz_n^{(k)}(a_n|h_n)\sz_n(b_n|h_n) \int_0^{s_{n+1}}D(\d t | i_n,a_n,b_n) \bigg) \Bigg\} \notag \\
& \qd + \Bigg\{ (1-\lz) q_k \one_{\{i\}}(i_0) \bigg[ \prod_{M=0}^{n-1} \hat{\pz}^{(k)}_M (a_M|h_M) \sz_{M} (b_M|h_M) \int_0^{s_{M+1}} Q( \d t, i_{M+1} | i_M,a_M,b_M) \bigg] \notag \\
& \qd \qd \times \bigg( \hat{\pz}_n^{(k)}(a_n|h_n)\sz_n(b_n|h_n) \int_0^{s_{n+1}}D(\d t | i_n, a_n, b_n) \bigg) \Bigg\} \notag\\
& = \one_{\{i\}} (i_0) \bigg[ \lz p_k \prod_{m=0}^n \pz_{m}^{(k)}(a_m|h_m)+ (1 - \lz) q_k \prod_{M=0}^n \hat{\pz}_M^{(k)} (a_M|h_M) \bigg] \notag \\
& \qd \times \bigg[ \prod_{m=0}^{n-1} \sz_{m} (b_m | h_m) \int_0^{s_{m+1}} Q( \d t, i_{m+1} |i_m, a_m, b_m) \bigg] \bigg( \sz_n(b_n | h_n) \int_0^{s_{n+1}} D(\d t | i_n, a_n, b_n) \bigg) \notag \\
& = \one_{\{i\}} (i_0) C_{n} (k, h_{n}, a_{n}) \bigg[ \prod_{m=0}^{n-1} \sz_{m} (b_m | h_m) \int_0^{s_{m+1}} Q( \d t, i_{m+1} |i_m, a_m, b_m) \bigg] \notag \\
& \qd \times \bigg( \sz_n(b_n | h_n) \int_0^{s_{n+1}} D(\d t |i_n, a_n, b_n) \bigg).
\end{align}
On the other hand, according to the definition of $\pz^\lz$ given in (\ref{s3-thm5-2}), we have
\begin{align}\lb{s3-thm5-4}
& \PP_{\lz p + (1- \lz) q}^{\pz^\lz, \sz} \lt( E_n \rt) \notag \\
& =  \one_{\{i\}}(i_0) (\lz p_k+(1-\lz)q_k) \bigg[ \prod_{M=0}^{n-1} \sz_{M}(b_M|h_M) \int_0^{s_{M+1}} Q( \d t, i_{M+1} |i_M, a_M, b_M) \bigg] \notag \\
& \qd \times \bigg( \sz_n(b_n | h_n) \int_0^{s_{n+1}} D(\d t | i_n, a_n, b_n) \bigg) \bigg[ \prod_{m=0}^n \pz_m^{\lz, (k)} (a_m | h_m) \bigg] \notag \\
& = \one_{\{i\}} (i_0) C_{n} (k, h_{n}, a_{n}) \bigg[ \prod_{m=0}^{n-1} \sz_{m} (b_m | h_m) \int_0^{s_{m+1}} Q( \d t, i_{m+1} |i_m, a_m, b_m) \bigg] \notag \\
& \qd \times \bigg( \sz_n(b_n | h_n) \int_0^{s_{n+1}} D(\d t |i_n, a_n, b_n) \bigg).
\end{align}
Combining with (\ref{s3-thm5-3}) and (\ref{s3-thm5-4}), we obtain that (\ref{s3-thm5-0}) holds for all measurable cylinder subsets $E_n$ of $(\ooz, \fq)$. Hence, (\ref{s3-thm5-0}) also holds for all $E \in \fq$.

\underline{Step 2, proof of (\ref{s3-thm5-0-0}).} Similarly, let $\{h_n \in \Ho_{n}, n\geq 0\}$ be the history of $G(p)$, and we define a sequence $\{D_n (h_n,b),n\geq 0, b\in \Bo\}$ as follow:
$$
D_n (h_n, b) = \lz \prod_{m=0}^n \sz_m (b| h_m) + (1- \lz) \prod_{m=0}^n \hat{\sz}_m (b | h_m).
$$
For convenience, denote $D_{-1} (h_{-1}, b) = 1$. Then, we define the policy $\sz^\lz = \lt\{ \sz^\lz_n, n \geq 0 \rt\}$ as
$$
\sz^\lz_n (b | h_n) :=
\left\{
  \begin{array}{ll}
    D_n (h_n, b) / D_{n-1} (h_{n-1}, b_{n-1}), & \hbox{$D_{n-1} (h_{n-1}, b_{n-1}) > 0$;} \\
    |B|^{-1}, & \hbox{otherwise.}
  \end{array}
\right.
$$
The rest only needs to verify that $\sz^\lz \in \ssz$ and (\ref{s3-thm5-0-0}). These are similar to the proof of Step 1, so the detailed calculations are ignored.
\deprf

\section{Existence of the value function and optimality equation}

The main conclusions of this section are the existence of the value function and the optimality equation. Under the Assumption \ref{ass1}, for each $p \in \pq (\Ko)$ and $i \in \So$, the value function $V^* (p, i)$ of $G(p)$ always exists, see Theorem \ref{s3-thm3} below. The other important conclusion is the optimality equation corresponding to the value function $V^* (p, i)$, see Theorem \ref{OEthm} below. For the convenience of describing the optimality equation, we introduce some notations.
\begin{itemize}
\item[\rm i)] For any $\gz = \{ \gz^{(k)} (\cdot) \in \pq (\Ao),  k \in \Ko \} \in \pq(A|K)$ and $a \in \Ao$, we define a mapping $\llz_{\gz, a}: \pq (\Ko) \to \pq (\Ko)$ as
\be\lb{llz}
\llz_{\gz, a} (p) (k) := \frac{\gz^{(k)} (a) p_{k}}{\sum_{l \in \Ko} \gz^{(l)} (a) p_l}, \qd \forall \ p \in \pq (\Ko), \ k \in K,
\de
in other words, $\llz_{\gz, a} (p) \in \pq(K)$.
\item[\rm ii)] Given $n \geq 0$, $\pz \in \ppz$ and $\sz \in \ssz$, we define the policies up to the $n$-th decision epoch as
\be\lb{nth-policy}
\pz|_n := \{\pz_m^{(k)}, k \in \Ko, 0 \leq m \leq n\}, \qd \sz|_n := \{ \sz_m, 0 \leq m \leq n\}.
\de
Denote by $\ppz[n]$ and $\ssz[n]$ the sets of all policies up to the $n$-th decision epoch with the form (\ref{nth-policy}), respectively. For the special case $n=0$, we have
\begin{align*}
\ppz[0] &= \lt\{ \mz | \mz^{(k)} (\cdot | i) \in \pq(\Ao), k \in \Ko, i \in \So \rt\}=\pq(A|K\times S);\\
\ssz[0] &= \lt\{ \nz | \nz (\cdot | i) \in \pq(\Bo), \text{ for all }i \in \So \rt\}=\pq(B|S).
\end{align*}
Similarly, we can define the sets $\ppz^D[n]$ and $\ssz^D[n]$ for the sets of deterministic policies up to the $n$-th decision epoch.
\item [\rm iii)] Let $\mathbb{M}$ be the set of real-valued functions $u$ defined on $\pq(K) \times S$. For each $\mz \in \ppz[0]$ and $\nz \in \ssz[0]$, we define operators  $T^{\mu, \nu} u (p, i)$, $\overline{T} u (p, i)$ and $\underline{T} u (p, i)$  from $\mathbb{M}$ to $\mathbb{M}$ as following:
\begin{align}
T^{\mu, \nu} u (p, i) & = \sum_{k \in \Ko} \sum_{a \in \Ao, b \in \Bo} p_{k} \mu^{(k)} (a | i) \nu (b | i) c (k, i, a, b) \int_0^{\infty} e^{-\az t} (1 - D(t| i, a, b)) \d t \notag \\
& \qd + \sum_{k \in \Ko} \sum_{a \in \Ao, b \in \Bo} p_{k} \mu^{(k)} (a | i) \nu (b | i) \sum_{j \in \So} \int_0^{\infty} e^{-\az t} Q (\d t, j | i, a, b) u (\llz_{\mz (\cdot| i), a} (p), j), \lb{T} \\
\overline{T} u (p, i) &= \inf_{\nu \in \ssz[0]} \sup_{\mu \in \ppz[0]} T^{\mu, \nu} u (p, i), \qqd \underline{T} u (p, i) = \sup_{\mu \in \ppz[0]} \inf_{\nu \in \ssz[0]} T^{\mu, \nu} u (p, i). \lb{TT}
\end{align}
\end{itemize}
\begin{thm}\lb{OEthm}
	{
		Suppose that Assumption \ref{ass1} holds. The value function of the semi-Markov game with incomplete information $G(p)$ satisfies the following optimality equation:
		$$
		V^* = \overline{T} V^* = \underline{T} V^*.
		$$
	}
\end{thm}

The proof of Theorem \ref{OEthm} is complicated and arranged in the last part of this section. It may be easy to read by showing how to find out the optimality equation step by step. At each step, we have either a proposition or a lemma. If one is in hurry, who may jump from here to the existences and the iterative algorithms of the optimal policies for \pla 1 (Theorem \ref{s4-thm1}, Algorithm \ref{algo-p1}) and the optimal policies for \pla 2 (Theorem \ref{s5-thm8}, Algorithm \ref{algo-p2}).

In order to prove the existence of value function $V^* (p, i)$ of $G(p)$, the key bridge is the expected discount reward up to the $n$-th decision epoch, which is defined as
\be\lb{Vn}
V_n (p, i, \pz, \sz) = \sum_{m=0}^n \EE_{p, i}^{\pz, \sz} \lt[ \frac{1}{\az} \lt(e^{-\az T_m} - e^{- \az T_{m+1}} \rt) c(\kz, X_m, A_m, B_m) \rt].
\de
The relationship between $V_n$ and $V$ is given below.
\begin{lem}\lb{s3-lem1}
	Suppose that Assumption \ref{ass1} holds. For each $\vz_0 >0$ there exists $N (\vz_0)$, which is independent of $p \in \pq (\Ko)$, $i \in \So$ and $(\pz, \sz) \in \ppz \times \ssz$, such that for all $n > N (\vz_0)$ it holds that
	$$
	V(p, i, \pz, \sz) - V_n (p, i, \pz, \sz) < \vz_0.
	$$
\end{lem}
\prf Since the payoff rate function $c(k, i, a, b)$ is nonnegative, $V_n(p,i,\pz,\sz)$ is non-decreasing with respect to $n \geq 0$. Moreover, the finiteness of $K$, $S$, $A$ and $B$ ensures that $c(k, i, a, b)$ is bounded, i.e. $c^* < \infty$. Hence, the monotone convergence theorem implies that
\begin{align*}
V (p, i, \pz, \sz) &= \EE_{p, i}^{\pz, \sz} \lt[\sum_{m=0}^{\infty} \int_{T_m}^{T_{m+1}} e^{-\az t} c(\kz, X_m, A_m, B_m) \d t \rt] \\
& = \sum_{m=0}^{\infty} \EE_{p, i}^{\pz, \sz} \lt[ \frac{1}{\az} \lt( e^{- \az T_m} - e^{- \az T_{m+1}} \rt) c(\kz, X_m, A_m, B_m) \rt].
\end{align*}
Hence, by the Proposition \ref{s2-prop1} we have
\begin{align}\lb{s3-lem1-1}
V (p, i, \pz, \sz) - V_n (p, i, \pz, \sz) &= \sum_{m = n+1}^\infty \EE_{p, i}^{\pz, \sz} \lt[ \frac{1}{\az} \lt(e^{-\az T_m} - e^{- \az T_{m+1}} \rt) c(\kz, X_m, A_m, B_m) \rt] \notag \\
&\leq \frac{c^*}{\az} \EE_{p, i}^{\pz, \sz} \lt[e^{-\az T_{n+1}} - e^{-\az T_{\infty}} \rt] = \frac{c^*}{\az} \EE_{p, i}^{\pz, \sz} \lt[e^{-\az T_{n+1}}\rt].
\end{align}
Using the same method given in (\ref{s2-prop2-2}), we obtain $\EE_{p, i}^{\pz, \sz} \lt[e^{-\az T_{n+1}}\rt] \leq \lt(1-\vz+\vz e^{-\az\delta} \rt)^{n+1}$, where $0< \vz < 1$ is introduced in Assumption \ref{ass1} and $\dz$ is a arbitrary positive constant. Note that $\bz:=1-\vz+\vz e^{-\az\delta} < 1$, then combining with (\ref{s3-lem1-1}), we obtain $N (\vz_0) = \lt| \log \frac{\az \vz_0}{c^*} \rt| \Big/ \lt| \log \bz \rt| + 1$ such that for all $n > N (\vz_0)$ is holds that
$$
V (p, i, \pz, \sz) - V_n (p, i, \pz, \sz) \leq \frac{c^*}{\az} \bz^{n+1} < \vz_0.
$$
The proof of this lemma is completed.
\deprf

Denote by $\overline{V_n} (p,i)$ and $\underline{V_n} (p,i)$ the upper and lower value of $V_n (p,i,\pz,\sz)$, respectively, i.e.,
$$
\overline{V_n} (p,i);=\inf_{\sigma\in \Sigma}\sup_{\pi\in \Pi}V_n(p,i,\pi,\sigma),\qd \underline{V_n} (p,i);=\sup_{\pi\in \Pi}\inf_{\sigma\in \Sigma}V_n(p,i,\pi,\sigma).
$$
When $\overline{V_n} = \underline{V_n}$, we say that the value function at the $n$-th decision epoch exists, which is denoted by $V^*_n$. The next lemma states that the value function $V_n^* (p, i)$ always exists. This gives us the direction to study the existence of the value function $V^* (p, i)$ of $G(p)$.

\begin{lem}\lb{s3-lem2}
	For each $n \geq 0$, we have
	\be\lb{s3-lem2-0}
	\overline{V_n} (p, i) = \underline{V_n} (p, i), \qd \forall p \in \pq (\Ko), \ i \in \So,
	\de
which means that the value function at the $n$-th decision epoch exists, denoted by $V_n^* (p, i)$.
\end{lem}

\prf The case $n=0$ is trivial. We consider the case $n=1$, the other cases can be analyzed and discussed in the same way. For each $\hat{\pz} \in \ppz[1]$ and $\hat{\sz} \in \ssz[1]$, we define
\begin{align*}
\pz &= \big\{ \pz_n^{(k)}, k \in K, n \geq 0 : \text{$\pz_0^{(k)}= \hat{\pz}_0^{(k)}$, $\pz_1^{(k)}= \hat{\pz}_1^{(k)}$ and $\pz_n^{(k)} \equiv |A|^{-1}$ for all $n \geq 2$ } \big\} \in \ppz, \\
\sz &= \big\{\sz_n, n \geq 0: \text{ $\sz_0 = \hat{\sz}_0$, $\sz_1 = \hat{\sz}_1$ and $\sz_n \equiv |B|^{-1}$ for all $n \geq 2$ }\big\} \in \ssz.
\end{align*}
Hence, we obtain $\pz|_1=\hat{\pz}$ and $\sz|_1=\hat{\sz}$. Then, for each $(\hat{\pz}, \hat{\sz}) \in \ppz[1] \times \ssz[1]$, define $\hat{V}_1 (p, i, \hat{\pz}, \hat{\sz}) := V_1 (p, i, \pz, \sz)$. Note a fact that the expected discount reward up to the $1$-st decision epoch $V_1$ only depends on the controlled processes $(T_m, X_m)$ and action processes $(A_m, B_m)$ with $m= 0, 1$ and is independent of the information of $m \geq 2$. Hence, we have
\begin{align*}
\sup_{\hat{\pz} \in \ppz[1]} \inf_{\hat{\sz} \in \ssz[1]} \hat{V}_1 (p, i, \hat{\pz}, \hat{\sz}) &= \sup_{\pz \in \ppz} \inf_{\sz \in \ssz} V_1 (p, i, \pz, \sz),\\
\inf_{\hat{\sz} \in \ssz[1]} \sup_{\hat{\pz} \in \ppz[1]} \hat{V}_1 (p, i, \hat{\pz}, \hat{\sz}) &= \inf_{\sz \in \ssz} \sup_{\pz \in \ppz} V_1 (p, i, \pz, \sz),
\end{align*}
which means that (\ref{s3-lem2-0}) holds for $n=1$ if and only if
\be\lb{s3-lem2-1}
\sup_{\hat{\pz} \in \ppz[1]} \inf_{\hat{\sz} \in \ssz[1]} \hat{V}_1 (p, i, \hat{\pz}, \hat{\sz})  = \inf_{\hat{\sz} \in \ssz[1]} \sup_{\hat{\pz} \in \ppz[1]} \hat{V}_1 (p, i, \hat{\pz}, \hat{\sz}).
\de

In the next step, we analysize the probability measure space on the deterministic policies, i.e. $\pq( \ppz^D[1])$ and $\pq (\ssz^D[1])$. The policies $\varphi \in \ppz^D[1]$ and $\psi \in \ssz^D [1]$ are described by the mappings $\varphi_n^{k}: \Ho_n \to \Ao$ and $\psi_n: \Ho_n \to \Bo$, where $k \in \Ko$ and $n=0, 1$. Since $\Ko$, $\Ao$, $\Bo$ and $\Ho_n$ are finite, the sets $\ppz^D[1]$ and $\ssz^D[1]$ are finite too. For any $x \in \pq (\ppz^D[1])$ and $y \in \pq (\ssz^D[1])$, we define
\begin{align*}
\widetilde{V}_1 (p, i, x, \hat{\sz}) &:= \sum_{\varphi \in \ppz^D[1]} x(\varphi) \hat{V}_1 (p, i, \varphi, \hat{\sz}), \qd \forall \hat{\sz} \in \ssz[1]; \qd \\
\widetilde{V}_1 (p, i, \hat{\pz}, y) &:= \sum_{\psi \in \ssz^D[1]} y(\psi) \hat{V}_1 (p, i, \hat{\pz}, \psi), \qd \forall \hat{\pz} \in \ppz[1]; \\
\widetilde{V}_1 (p, i, x, y) &:= \sum_{\varphi \in \ppz^D[1]} \sum_{\psi \in \ssz^D[1]} x(\varphi) y(\psi) \hat{V}_1 (p, i, \varphi, \psi).
\end{align*}
By the definition, we obtain that $\widetilde{V}_1 (p, i, \cdot, \cdot)$ is a bilinear function on $\pq (\ppz^D[1]) \times \pq (\ssz^D[1])$.

The relationship between the randomized history-dependent policies in $\ppz[1]$ (resp. $\ssz[1]$) and the probability measures of the deterministic policies in $\pq( \ppz^D[1])$ (resp. $\pq (\ssz^D[1])$) is given by two steps. Firstly, we construct an equivalent randomized history-dependent policy from each $x \in \pq (\ppz^D[1])$. In detail, the policy $\pz^{x} = \{ {\pz}_0^{{x}, (k)},  {\pz}_1^{{x}, (k)}, k \in \Ko \}$ induced by $x \in \pq (\ppz^D[1])$ is defined as
\begin{align*}
{\pz}_0^{{x}, (k)}(a | h_0) &= \sum_{\varphi \in \ppz^D[1]} x(\varphi) \one_{\{ \varphi_0^{(k)} (h_0)\}} (a), \\
{\pz}_1^{{x}, (k)} (a | h_1) & =
\left\{
\begin{array}{ll}
\dps\frac{1}{{\pz}_0^{{x}, (k)} (a_0 | h_0)} {\dps\sum_{\varphi \in \ppz^D [1]} x(\varphi) \one_{\{\varphi_0^{(k)} (h_0)\}} (a_0) \one_{\{ \varphi_1^{(k)} (h_1)\}} (a)}, & \hbox{${\pz}_0^{ {x}(k)} (a_0 | h_0) > 0$,} \\
|A|^{-1}, & \hbox{otherwise,}
\end{array}
\right.
\end{align*}
Using the result of \cite[Theorem D.1]{SS02}, we obtain $\pz^{x} \in \ppz[1]$ and
\be\lb{s3-lem2-2-1}
\hat{V}_1 (p, i, \pz^{x}, \hat{\sz})=\widetilde{V}_1 (p, i, x, \hat{\sz})  ,\qd  \forall {\,}\hat{\sz} \in \ssz[1].
\de
According to the similar construction, given any $y \in \pq (\ssz^D[1])$, there exists ${\sz^{y}} \in \ssz[1]$ such that
\be\lb{s3-lem2-2-2}
\hat{V}_1 (p, i, \hat{\pz}, \sz^{y}) = \widetilde{V}_1 (p, i, \hat{\pz}, y),\qd \forall{\,} \hat{\pz}\in \ppz[1].
\de
Secondly, we construct an equivalent measure in $\pq (\ppz^D[1])$ from each randomized history-dependent policy. Given any $\hat{\pz} = \{ \hat{\pz}_0^{(k)}, \hat{\pz}_1^{(k)}, k \in \Ko\} \in \ppz[1]$, we define
$$
x^{\hat{\pz}} (\varphi) := \prod_{k \in \Ko, \ i_0 \in \Ho_0} \hat{\pz}_0^{(k)} \lt(\varphi_0^{(k)} (i_0) \big| i_0 \rt) \times \prod_{l \in K, \ h_1 \in H_1} \hat{\pz}_1^{(l)} \lt( \varphi_1^{(l)} (h_1) \big| h_1 \rt), \qd  \forall \,\varphi \in \ppz^D[1].
$$
Using the result of \cite[Theorem D.1]{SS02} again, we have $x^{\hat{\pz}}\in \pq (\ppz^D[1])$ and
\begin{align}\lb{s3-lem2-3-1}
\widetilde{V}_1(p, i, x^{\hat{\pz}}, \hat{\sz})=\hat{V}_1(p, i, \hat{\pz}, \hat{\sz}) ,\qd \forall\, \hat{\sz} \in \ssz[1].
\end{align}
Similarly, given any $\hat{\sz} \in \ssz[1]$, there exists $y^{\hat{\sz}} \in \pq (\ssz^D[1])$ such that
\be\lb{s3-lem2-3-2}
\widetilde{V}_1 (p, i, \hat{\pz}, y^{\hat{\sz}})=\hat{V}_1(p, i, \hat{\pz}, \hat{\sz}) , \qd \forall \, \hat{\pz} \in \ppz[1].
\de
Hence, the formulas (\ref{s3-lem2-2-1}) - (\ref{s3-lem2-3-2}) and the definition of $\widetilde{V}_1$ imply that
\begin{align*}
\sup_{\hat{\pz} \in \ppz[1]} \inf_{\hat{\sz} \in \ssz[1]} \hat{V}_1 (p, i, \hat{\pz}, \hat{\sz}) &= \sup_{x \in \pq (\ppz^D[1])} \inf_{y \in \pq (\ssz^D[1])} \widetilde{V}_1 (p, i, x, y), \\
\inf_{\hat{\sz} \in \ssz[1]} \sup_{\hat{\pz} \in \ppz[1]} \hat{V}_1 (p, i, \hat{\pz}, \hat{\sz}) &= \inf_{y \in \pq (\ssz^D[1])} \sup_{x \in \pq (\ppz^D[1])} \widetilde{V}_1 (p, i, x, y).
\end{align*}
Note that $\ppz^D[1]$ (resp. $\ssz^D[1]$) is finite set, then $\pq (\ppz^D[1])$ (resp. $\pq (\ssz^D[1])$) is convex and compact set, see \cite{B68}. Moreover, since $\widetilde{V}_1(p, i, \cdot, \cdot)$ is a bilinear function with respect to $(x, y) \in \pq (\ppz^D[1]) \times \pq (\ssz^D[1])$, the von-Neumann minimax theorem (see \cite[Theorem 1.2.3]{B13}) implies that there exists $(x^*, y^*) \in \pq (\ppz^D[1]) \times \pq (\ssz^D[1])$ such that
\begin{align*}
& \sup_{x \in \pq (\ppz^D[1])} \inf_{y \in \pq (\ssz^D[1])} \widetilde{V}_1 (p, i, x, y) = \sup_{x \in \pq (\ppz^D[1])} \widetilde{V}_1 (p, i, x, y^*) \\
& \qd = \inf_{y \in \pq (\ssz^D[1])} \widetilde{V}_1 (p, i, x^*, y) = \inf_{y \in \pq (\ssz^D[1])} \sup_{ x\in \pq (\ppz^D[1])}\widetilde{V}_1 (p, i, x, y).
\end{align*}
Hence, the above equation implies that (\ref{s3-lem2-1}) holds, and then (\ref{s3-lem2-0}) holds equivalently.
\deprf

Lemma \ref{s3-lem2} ensures the existence of the value function $V_n^* (p, i)$ at the $n$-th decision epoch. Intuitively, the value function $V^* (p, i)$ of $G(p)$ can be constructed by letting $n \to \infty$. Details are given below.

\begin{thm}\lb{s3-thm3}
	Suppose that Assumption \ref{ass1} holds. Fixed any $p \in \pq (\Ko)$, the value function of $G(p)$ exists and satisfies
	$$
	V^* (p, i) = \lim_{n \to \infty} V_n^* (p, i) , \qqd \forall \, i \in \So.
	$$
\end{thm}

\prf For each $n \geq 0$ and $(\pz, \sz) \in \ppz \times \ssz$, the fact that $V_n (p, i, \pz, \sz) \leq V(p, i, \pz, \sz)$ implies
\be\lb{s3-thm3-1}
\underline{V_n}(p,i)=\sup_{\pz \in \ppz} \inf_{\sz \in \ssz} V_n (p, i, \pz, \sz) \leq \sup_{\pz \in \ppz} \inf_{\sz \in \ssz} V(p, i, \pz, \sz)=\underline{V}(p,i), \qd \forall i \in \So.
\de
According to Lemma \ref{s3-lem1}, for each $\vz_0 > 0$, there exists $N (\vz_0)$ such that for all $n \geq N (\vz_0)$, we have $$
\overline{V}(p, i)= \inf_{\sz \in \ssz}\sup_{\pz \in \ppz} V(p, i, \pz, \sz) \leq\inf_{\sz \in \ssz} \sup_{\pz \in \ppz} ( V_n (p, i, \pz, \sz)+\vz_0) =\overline{V_n} (p, i) + \vz_0.
$$
 Combining with Lemma \ref{s3-lem2} and (\ref{s3-thm3-1}), we obtain
\be\lb{s3-thm3-2}
V^*_n (p, i) \leq \underline{V}(p, i) \leq \overline{V}(p, i ) \leq V^*_n (p, i) + \vz_0.
\de
Moreover, by the definition of $V_n^* (p, i)$, it is bounded and nondecreasing with respect to $n$. Hence, the limit of $\lim_{n \to \infty} V_n^* (p, i)$ exists, which is denoted by $V^* (p, i)$. Then, (\ref{s3-thm3-2}) implies that
$$
V^* (p, i) \leq \underline{V}(p, i) \leq \overline{V}(p, i ) \leq V^* (p, i) + \vz_0.
$$
By the arbitrariness of $\vz_0$, we obtain $\lim_{n \to \infty} V_n^* (p, i)=V^* (p, i) =\underline{V}(p, i) =\overline{V}(p, i )$. Hence, the existence of the value function of $G(p)$ is completed.
\deprf

Fix any $i \in \So$, $V^* (\cdot, i)$ and $V_n^*(\cdot, i)$ are functions defined on the space $(\pq (\Ko), \| \cdot \|)$. Naturally, we will give some properties of the value function $V^* (p, i)$ and $V^*_n(p,i)$ in the next step, such as continuity (see Lemma \ref{s3-lem4}) and concavity (see Lemma \ref{s3-cor6}).

\begin{lem}\lb{s3-lem4}
For each $i \in \So$ and $n \geq 0$, the value function at the $n$-th decision epoch $V^*_n (\cdot, i)$ is Lipschitz continuous with respect to $p \in \pq (\Ko)$, i.e., there exists a constant $C > 0$ such that
$$
|V^*_n(p, i) - V^*_n(q, i)| \leq C \| p - q \|, \qd \forall p, \ q \in \pq (\Ko).
$$
Moreover, under Assumption \ref{ass1}, the value function $V^* (\cdot, i)$ is also Lipschitz continuous with respect to $p \in \pq (\Ko)$.	
\end{lem}

\prf By the linearity property given in (\ref{linear}), we obtain that $V_n(p, i, \pz, \sz) = \sum_{k \in K} p_k V_n(\dz_k, i, \pz, \sz)$, where $\dz_k$ is the Dirac measure on $\pq (\Ko)$. For any $p, q \in \pq (\Ko)$ and $(\pi,\sz)\in \Pi\times \Sigma$, we have
\begin{align}
& |V_n(p, i, \pz, \sz) - V_n(q, i, \pz, \sz)| = \lt| \sum_{k \in \Ko} p_k V_n(\dz_k, i, \pz, \sz) - \sum_{k \in \Ko} q_k V_n(\dz_k, i, \pz, \sz) \rt| \notag \\
& \qd \leq \sum_{k \in K} |p_k - q_k| \EE_{\dz_k, i}^{\pz, \sz} \lt[\int_0^{T_{n+1}}e^{-\az t} c(\kz, X_t, A_t, B_t) \d t \rt] \leq \frac{c^*}{\az} \| p - q \|. \lb{s3-lem4-1}
\end{align}
According to Lemma \ref{s3-lem2}, the value functions $V^*_n(p, i)$ and $V^*_n(q, i)$ exist. Since the bound of (\ref{s3-lem4-1}) is independent of policies $(\pz, \sz)$, we obtain $|V^*_n(p, i) - V^*_n(q, i)| \leq \frac{c^*}{\az} \| p - q \|$, which means $V_n^*(\cdot, i)$ is Lipschitz continuous with Lipschitz constant $C= c^* / \az$. Moreover, noting that $C$ is independent of $i \in \So$ and $n$, under Assumption \ref{ass1}, Theorem \ref{s3-thm3} implies that
$$
|V^*(p,i) - V^*(q,i)| \leq C \|p-q\|,
$$
which means that $V^* (p, i)$ is also Lipschitz continuous with respect to $p \in \pq (K)$.
\deprf

\begin{lem}\lb{s3-cor6}
For each $i \in \So$ and $n \geq 0$, the value function at the $n$-th decision epoch $V^*_n (\cdot, i)$ is concave on $\pq (\Ko)$. Furthermore, under Assumption \ref{ass1}, the value function $V^* (\cdot, i)$ is also concave.
\end{lem}

\prf Fixed any $n \geq 0$ and $p, q \in \pq (\Ko)$, by the definition of $\underline{V_n}$ and Lemma \ref{s3-lem2} , for each $\vz_0 >0$, there exist $\pz \in \ppz$ and $\hat{\pz} \in \ppz$ such that
\begin{align*}
\inf_{\sz \in \ssz} V_n (p, i, \pz, \sz) \geq \underline{V_n} (p, i) - \vz_0 = V^*_n(p,i) - \vz_0, \\
\inf_{\sz \in \ssz} V_n (q, i, \hat{\pz}, \sz) \geq \underline{V_n} (q, i) - \vz_0 =V^*_n(q,i) - \vz_0.
\end{align*}
For each $\lz \in [0, 1]$, Proposition \ref{s3-thm5} implies that there exists ${\pz^\lz} \in \ppz$ such that
$$
V_n (\lz p + (1- \lz) q, i, {\pz^\lz}, \sz) = \lz V_n (p, i, \pz, \sz) + (1- \lz) V_n (q, i, \hat{\pz}, \sz), \qd \forall \sz \in \ssz.
$$
Hence, we have
\begin{align*}
&V^*_n (\lz p + (1- \lz) q, i) \geq \inf_{\sz \in \ssz} V_n (\lz p + (1- \lz) q, {\pz^\lz}, \sz) \\
& \qd \geq \lz \inf_{\sz \in \ssz} V_n (p, i, \pz, \sz) + (1- \lz) \inf_{\sz \in \ssz} V_n (q, i, \hat{\pz}, \sz) \geq \lz {V_n^*} (p, i) + (1- \lz) {V_n^*} (q, i) - \vz.
\end{align*}
By the arbitrariness of $\vz$, we have $V^*_n (\lz p + (1- \lz) q, i) \geq \lz V^*_n (p, i) + (1- \lz) V^*_n (q, i)$, i.e,  $V^*_n  (\cdot, i)$ is a concave function on $\pq (\Ko)$. Furthermore, under Assumption \ref{ass1}, Theorem \ref{s3-thm3}  gives that $V^* (p, i) = \lim_{n \to \infty}V_n^* (p, i)$, which implies that $V^* (p, i)$ is also concave.
\deprf

Similar to the discrete time Markov decision process (see \cite{BR11} for instance), the value function of $G(p)$ may satisfy some recursive relations. However, in semi-Markov processes, this relationship is more complex. In the following result, we introduce the recursive relation of $V_n(p, i, \pz, \sz)$ and $V(p,i,\pz,\sz)$. Before stating the result, we give the definition of one-step forward policies. Fixed $(i,a,b)\in S\times A\times B$, the one-step forward policy of $\pz \in \ppz$ is denoted by $^{(i, a, b)}\pz=\{{^{(i, a, b)}\pz}_n^{(k)},n\geq 0, k\in K\}$, which satisfies
\be\lb{1-pz}
{^{(i, a, b)}\pz}_n^{(k)} (\cdot| h_n) := \pz_{n+1}^{(k)} (\cdot | i, a, b, h_{n}), \qd \forall h_n \in \Ho_n,
\de
where $(i,a,b, h_{n})=(i,a,b,i_0,a_0,b_0,i_1,\ldots,a_{n-1},b_{n-1},i_n)\in \Ho_{n+1}$. The one-step forward policy $^{(i, a, b)}\sz=\{{^{(i, a, b)}\sz}_n, n\geq 0\}$ of $\sz \in \ssz$ is defined in similar way.

\begin{lem}\lb{s3-lem7}
	For any $n\geq 0$, the function $V_n(p, i, \pz, \sz)$ satisfies
	\begin{align}\lb{s3-lem7-0-0}
	& V_{n+1} (p, i, \pz, \sz) = \sum_{k \in \Ko} \sum_{a \in \Ao, \ b \in \Bo} p_k \pz_0^{(k)} (a | i) \sz_0 (b | i) c (k, i, a, b) \int_0^\infty e^{-\az t} (1 - D(t | i, a, b)) \d t \notag \\
	& + \sum_{k \in \Ko} \sum_{a \in \Ao, \  b \in \Bo} p_k \pz_0^{(k)} (a|i) \sz_0 (b|i) \sum_{j \in \So} \int_0^\infty e^{-\az t} Q (\d t, j| i, a, b) V _n(\llz_{\pz_0 (\cdot | i), a} (p), j, {^{(i, a, b)}\pz}, {^{(i, a, b)}\sz}),
	\end{align}
	where $\llz_{\pz_0 (\cdot | i), a} (p)$ is the probability measure on $\pq (K)$ defined in (\ref{llz}).
Moreover, under Assumption \ref{ass1},  the function $V(p, i, \pz, \sz)$ satisfies
	\begin{align}\lb{s3-lem7-0}
	& V (p, i, \pz, \sz) = \sum_{k \in \Ko} \sum_{a \in \Ao, \ b \in \Bo} p_k \pz_0^{(k)} (a | i) \sz_0 (b | i) c (k, i, a, b) \int_0^\infty e^{-\az t} (1 - D(t | i, a, b)) \d t \notag \\
	& + \sum_{k \in \Ko} \sum_{a \in \Ao, \  b \in \Bo} p_k \pz_0^{(k)} (a|i) \sz_0 (b|i) \sum_{j \in \So} \int_0^\infty e^{-\az t} Q (\d t, j| i, a, b) V (\llz_{\pz_0 (\cdot | i), a} (p), j, {^{(i, a, b)}\pz}, {^{(i, a, b)}\sz}).
	\end{align}
\end{lem}

\prf Fix any $n\geq 0$. By (\ref{E[f]}) and (\ref{Vn}), we have that
\begin{align}\lb{s3-lem7-1}
V_{n+1} (p, i, \pz, \sz) & = \sum_{k \in \Ko} \sum_{a \in \Ao, \  b \in \Bo} p_k \pz_0^{(k)} (a|i) \sz_0 (b|i) c(k, i, a, b) \lt( \int_0^\infty e^{-\az t} \lt(1 - D (t | i, a, b) \rt) \d t \rt) \notag \\
& \qd + \EE_{p, i}^{\pz, \sz} \lt[ \sum_{m=1}^{n+1} \frac{1}{\az} (e^{-\az T_m} - e^{-\az T_{m+1}}) c(\kz, X_m, A_m, B_m) \rt].
\end{align}
To calculate the second item of (\ref{s3-lem7-1}), we give some preparations. Based on the Bayes formula, we have
\begin{align*}
& \PP_{p, i}^{\pz, \sz} (\kz = k |  A_0 = a, B_0 = b, X_1=j, T_1\in [0,s_1]) \notag \\
& \qd = \frac{p_k \pz_0^{(k)} (a|i) \sz_0 (b|i)\int_0^{s_1} Q(\d t,j|i,a,b)}{\sum_{l \in \Ko} p_l \pz_0^{(l)}(a|i) \sz_0(b|i)\int_0^{s_1} Q(\d t,j|i,a,b)} = \llz_{\pz_0 (\cdot | i), a} (p) (k),
\end{align*}
which implies that for any function $f:\Ko \to \RR_+$, it holds that
\be\lb{s3-lem7-3}
\EE^{\pz,\sz}_{p,i}[f(\kz)|A_0,B_0,X_1,T_1]=\sum_{k\in \Ko}\llz_{\pz_0 (\cdot | i), A_0} (p) (k)f(k).
\de
Using the expression (\ref{E[f]}) again, for any $1\leq m \leq n$, we calculate the conditional expectation as following
\begin{align}\lb{s3-lem7-4}
& \EE_{\dz_k, i}^{\pz, \sz} \lt[ \lt( e^{-\az (T_m - T_1)} - e^{-\az (T_{m+1} - T_1)} \rt) c(\kz, X_m, A_m, B_m) \big|A_0 =a, B_0= b, X_1 = j,\kz=k \rt] \notag \\
& = \sum_{i_0}\one_{\{j\}}(i_0)\sum_{a_0, b_0} \pz_1^{(k)} (a_0 | i,a,b,h_{0}) \sz_1(b_0 | i,a,b,h_{0}) \sum_{i_1} \int_0^\infty Q(\d t_1, i_1| i_0, a_0, b_0)\notag \\
& \qd \cdots \sum_{a_{m-1}, b_{m-1}} \pz_{m}^{(k)} (a_{m-1}| i,a,b,h_{m-1}) \sz_{m} (b_{m-1} | i,a,b,h_{m-1})  \notag \\
& \qd \sum_{i_{m}} \int_0^\infty Q(\d t_{m}, i_{m} | i_{m-1}, a_{m-1}, b_{m-1})c(k, i_{m-1}, a_{m-1}, b_{m-1}) \lt( e^{-\alpha \sum_{l=1}^{m-1} t_l} - e^{-\alpha \sum_{l=1}^{m} t_l} \rt) \notag \\
& = \sum_{i_0} \one_{\{j\}} (i_0) \sum_{a_0, b_0} {^{(i,a,b)} \pz_0^{(k)}} (a_0|h_0) {^{(i,a,b)}\sz_0} (b_0|h_0) \sum_{i_1} \int_0^\infty Q(\d t_1, i_1|i_0,a_0,b_0) \notag \\
& \qd \cdots \sum_{a_{m-1}, b_{m-1}} {^{(i,a,b)} \pz_{m-1}^{(k)}} (a_{m-1} | h_{m-1}) {^{(i,a,b)}\sz_{m-1}} (b_{m-1} | h_{m-1}) \notag \\
& \qd \sum_{i_{m}} \int_0^\infty Q(\d t_{m}, i_{m}|i_{m-1}, a_{m-1}, b_{m-1}) c(k, i_{m-1}, a_{m-1}, b_{m-1}) \lt(e^{-\alpha\sum_{l=1}^{m-1} t_l} - e^{-\alpha \sum_{l=1}^{m} t_l} \rt) \notag \\
& = \EE_{\dz_k, j}^{^{(i, a, b)}\pz, ^{(i, a, b)}\sz} \lt[ \lt( e^{-\az T_{m-1}} - e^{-\az T_{m}} \rt) c(\kz, X_{m-1}, A_{m-1}, B_{m-1}) \rt].
\end{align}
Hence, using (\ref{s3-lem7-3}) and (\ref{s3-lem7-4}), the second item of (\ref{s3-lem7-1}) satisfies
\begin{align}\lb{s3-lem7-5}
& \EE_{p, i}^{\pz, \sz} \Bigg[ \sum_{m=1}^{n+1} \frac{1}{\az} \lt(e^{-\az T_m} - e^{-\az T_{m+1}} \rt) c(\kz, X_m, A_m, B_m) \Bigg] \notag \\
=&\EE_{p, i}^{\pz, \sz} \Bigg[\EE_{p, i}^{\pz, \sz} \Bigg[ \sum_{m=1}^{n+1} \frac{1}{\az} \lt(e^{-\az T_m} - e^{-\az T_{m+1}} \rt) c(\kz, X_m, A_m, B_m) \Bigg|A_0,B_0,\kz,X_1,T_1\Bigg]\Bigg] \notag\\
=&\EE_{p, i}^{\pz, \sz} \Bigg[ e^{-\az T_1} \sum_{m=1}^{n+1} \EE_{p, i}^{\pz, \sz}\Bigg[  \frac{1}{\az} \lt(e^{-\az (T_m-T_1)} - e^{-\az (T_{m+1}-T_1)} \rt) c(\kz, X_m, A_m, B_m) \Bigg|A_0,B_0,\kz,X_1,T_1\Bigg]\Bigg] \notag\\
=&\EE_{p, i}^{\pz, \sz} \bigg[e^{-\az T_1}\EE_{p, i}^{\pz, \sz} \bigg[V_n(\dz_{\kz},X_1, {^{(i,A_0,B_0)}\pz},{^{(i,A_0,B_0)}\sz})\bigg|A_0,B_0,X_1,T_1\bigg]\bigg]\notag\\
=&\EE_{p, i}^{\pz, \sz}\bigg[ e^{-\az T_1} \sum_{k\in K}\llz_{\pz_0 (\cdot | i), A_0} (p) (k)V_n(\dz_{k},X_1, {^{(i,A_0,B_0)}\pz},{^{(i,A_0,B_0)}\sz})\bigg] \notag\\
= &\sum_{k} \sum_{a, b} p_k \pz_0^{(k)} (a | i) \sz_0 (b | i) \sum_{j} \int_0^\infty e^{-\az t} Q (\d t, j | i, a, b) V_n (\llz_{\pz_0 (\cdot | i), a} (p), j, {^{(i, a, b)}\pz}, {^{(i, a, b)}\sz}).
\end{align}
Hence, we complete the formula (\ref{s3-lem7-0-0}) by using (\ref{s3-lem7-1}) and (\ref{s3-lem7-5}). Finally, under Assumption \ref{ass1}, Theorem \ref{s3-thm3} and monotone convergence theorem give the formula (\ref{s3-lem7-0}) by letting $n \to \infty $ in the formula (\ref{s3-lem7-0-0}).
\deprf

Having these preparations (Theorem \ref{s3-thm3}, Lemma \ref{s3-lem4}, Corollary \ref{s3-cor6} and Lemma \ref{s3-lem7}) at hand, we can give the optimality equation corresponding to the value function $V^* (p, i)$, which is the main work of this section.

{\bf Proof of Theorem \ref{OEthm}.}
\underline{Step 1.} According to Theorem \ref{s3-thm3}, the value function $V^*(p, i)$ exists for each $p \in \pq (\Ko)$ and $i \in \So$. We now show that $\overline{T} V^* = \underline{T} V^*$, i.e.
\be\lb{s3-OE-1}
\inf_{\nu \in \ssz[0]} \sup_{\mu \in \ppz[0]} T^{\mu, \nu} V^* (p, i) = \sup_{\mu \in \ppz[0]} \inf_{\nu \in \ssz[0]} T^{\mu, \nu} V^* (p, i), \qd\forall \,p \in \pq (\Ko), i \in \So.
\de
Fix any $r \in \RR$, $\mu_0 \in \ppz[0]$ and $\nu_0 \in \ssz[0]$, define sets as
\begin{align*}
D(\nu_0) &:= \lt\{\mu \in \ppz[0] : T^{\mu, \nu_0} V^* (p, i) \geq r \rt\} \subseteq \ppz[0], \\
E (\mu_0) &:=  \lt\{\nu \in \ssz[0] : T^{\mu_0, \nu} V^* (p, i) \leq r \rt\} \subseteq \ssz[0].
\end{align*}
When the sets $D(\nu_0)$ and $E (\mu_0)$ are the convex closed subsets of $\ppz[0]$ and $\ssz[0]$ respectively, we can use the Sion minimax theorem (see \cite[Theorem A.7]{SS02}) immediately to show the equation(\ref{s3-OE-1}).

Firstly, we consider the $E (\mu_0)$. Given any $\{ \nu_n \}_{n \geq 0} \subseteq E (\mu_0)$ satisfying $\nu_n \to \nu$ as $n \to \infty$. Since $\Ko$, $\Ao$ and $\Bo$ are finite sets, the limit and the summations in $T^{\mz, \nz}$ (see the definition in (\ref{T})) can be exchanged, then we have $\lim_{n \to \infty} T^{\mu_0, \nu_n} V^* (p, i) = T^{\mz_0, \nu} V^* (p, i)$, which means that $\nu \in E (\mu_0)$. By the definition in (\ref{T}), the operator $T^{\mu_0, \nu}$ is linear with respect to $\nu \in \ssz[0]$. Then, for each $\nu, \hat{\nu} \in E (\mu_0)$ and $\lz \in [0, 1]$, it holds that
$$
T^{\nu_0, \lz \nu+ (1- \lz) \hat{\nu}} V^* (p, i) = \lz T^{\mz_0, \nu} V^*(p, i) + (1- \lz) T^{\mz_0, \hat{\nu}} V^*(p, i),
$$
which implies $\lz \nu+ (1- \lz) \hat{\nu} \in E (\mu_0)$. Hence, $E (\mu_0)$ is a convex closed subsets of $\ssz[0]$.

Secondly, we consider the $D(\nu_0)$. Given any $\{ \mu_n \}_{n \geq 0} \subseteq D (\nu_0)$ satisfying $\mu_n \to \mu$ as $n \to \infty$, we have
$$
\lim_{n \to \infty} \llz_{\mu_n (\cdot | i), a} (p) (k) = \dps\frac{p_k \lim_{n \to \infty} \mu_n^{(k)} (a | i)}{\sum_{l \in \Ko} p_l \lim_{n \to \infty} \mu_n^{(l)} (a|i)} = \llz_{\mu(\cdot | i), a} (p) (k).
$$
Combining with the continuity of $V^* (p, i)$ given in Lemma \ref{s3-lem4}, we have $T^{\mu, \nu_0} V^*(p, i) = \lim_{n \to \infty} T^{\mu_n, \nu_0} V^* (p, i)$, which means $D(\nu_0)$ is a closed set. Next, for each $\mu, \hat{\mu} \in D (\nu_0)$ and $\lz \in [0, 1]$, we show that
\be\lb{s3-OE-2}
T^{\lz \mu + (1 - \lz) \hat{\mu}, \nu_0} V^* (p, i) \geq \lz T^{\mu, \nu_0} V^* (p, i) + (1- \lz) T^{\hat{\mu}, \nu_0} V^* (p, i).
\de
To do so, fix any $i \in \So$, denote by $M_k^{a} = p_k \mu^{(k)} (a | i)$ and $\hat{M}_k^{a} = p_k \hat{\mu}^{(k)} (a | i)$. Note that
\begin{align*}
& \llz_{\lz \mu (\cdot | i) + (1- \lz) \hat{\mu} (\cdot | i), a} (p) (k) \\
=& \frac{\sum_{l \in \Ko}  M_l^{a}}{\sum_{l \in \Ko} \lt[ \lz M_l^{a} + (1- \lz) \hat{M}_l^{a} \rt]} \frac{ \lz M_k^{a}}{\sum_{l \in \Ko} M_l^a} + \frac{\sum_{l \in \Ko}  \hat{M}_l^{a}}{\sum_{l \in \Ko} \lt[ \lz M_l^{a} + (1- \lz) \hat{M}_l^{a} \rt]} \frac{(1- \lz)\hat{M}_k^{a}}{\sum_{l \in \Ko} \hat{M}_l^{a}} \\
=& \frac{\sum_{l \in \Ko} \lz M_l^{a}}{\sum_{l \in \Ko} \lt[ \lz M_l^{a} + (1- \lz) \hat{M}_l^{a} \rt]} \llz_{\mu(\cdot | i), a} (p) (k) + \frac{\sum_{l \in \Ko} (1- \lz) \hat{M}_l^{a}}{\sum_{l \in \Ko} \lt[ \lz M_l^{a} + (1- \lz) \hat{M}_l^{a} \rt]} \llz_{\hat{\mu}(\cdot | i), a} (p) (k).
\end{align*}
By Lemma \ref{s3-cor6}, $V^* (p, i)$ is a concave function on $\pq (\Ko)$, we have
\begin{align*}
V^* (\llz_{\lz \mu(\cdot | i) + (1- \lz) \hat{\mu} (\cdot | i), a} (p), i) \geq& \frac{\sum_{l \in \Ko} \lz M_l^{a}}{\sum_{l \in \Ko} \lt[ \lz M_l^{a} + (1- \lz) \hat{M}_l^{a} \rt]} V^* (\llz_{\mu(\cdot | i), a} (p), i) \\ &+ \frac{\sum_{l \in \Ko} (1- \lz) \hat{M}_l^{a}}{\sum_{l \in \Ko} \lt[ \lz M_l^{a} + (1- \lz) \hat{M}_l^{a} \rt]} V^* (\llz_{\hat{\mu}(\cdot | i), a} (p) , i).
\end{align*}
Hence, we obtain
\begin{align*}
& T^{\lz \mu + (1 - \lz) \hat{\mu}, \nu_0} V^* (p, i) \\
=& \sum_{k} \sum_{a, b} \lt( \lz M_k^{a} + (1 - \lz) \hat{M}_k^{a} \rt) \nu_0 (b | i) c (k, i, a, b) \int_0^\infty e^{-\az t} (1- D(t|i, a, b)) \d t \\
&  + \sum_{k} \sum_{a, b} \lt( \lz M_k^{a} + (1- \lz) \hat{M}_k^{a} \rt) \nu_0 (b | i) \sum_{j}  \int_0^\infty e^{-\az t} Q (\d t, j | i, a, b) V^* (\llz_{\lz \mu (\cdot | i) + (1- \lz) \hat{\mu} (\cdot|i), a} (p), j)  \\
\geq &\sum_{k} \sum_{a, b} \lt( \lz M_k^{a} + (1 - \lz) \hat{M}_k^{a} \rt) \nu_0 (b | i) c (k, i, a, b) \int_0^\infty e^{-\az t} (1- D(t|i, a, b)) \d t \\
&  + \sum_{a, b} \nu_0 (b | i) \sum_{j} \int_0^\infty e^{-\az t} Q (\d t, j | i, a, b) \Bigg[ \lz \bigg(\sum_{k} M_k^{a} \bigg) V^* (\llz_{\mu (\cdot | i) , a} (p), j) \\
&  + (1- \lz) \bigg(\sum_{k} \hat{M}_k^{a} \bigg) V^* (\llz_{\hat{\mu} (\cdot | i) , a} (p), j) \Bigg] \\
=& \lz T^{\mu, \nu_0} V^* (p, i) + (1 - \lz) T^{\hat{\mu}, \nu_0} V^* (p, i).
\end{align*}
Hence, (\ref{s3-OE-2}) holds and $D (\nu_0)$ is a convex closed subsets of $\ppz[0]$.

Finally, using the Sion minimax theorem in \cite[Theorem A.7]{SS02}, we obtain that the formula (\ref{s3-OE-1}) holds, i.e. $\overline{T} V^* = \underline{T} V^*$. Moreover, under the compactness of $\ppz[0]$ and $\ssz[0]$, the Sion minimax theorem in \cite[Theorem A.7]{SS02} also ensures that there exist $\mu_{p} \in \ppz[0]$ and $\nu_{p} \in \ssz[0]$ such that
\begin{align}
\underline{T} V^* (p, i) &= \sup_{\mu \in \ppz[0]} \inf_{\nu \in \ssz[0]} T^{\mu, \nu} V^* (p, i) = \inf_{\nu \in \ssz[0]} T^{\mu_{p}, \nu} V^* (p, i), \qd \forall i\in S; \lb{s3-OE-3} \\
\overline{T} V^* (p, i) &= \inf_{\nu \in \ssz[0]} \sup_{\mu \in \ppz[0]} T^{\mu, \nu} V^* (p, i) = \sup_{\mu \in \ppz[0]} T^{\mu, \nu_{p}} V^* (p, i) \qd \forall i\in S. \lb{s3-OE-4}
\end{align}

\underline{Step 2.} Next, we will show that
\be\lb{s3-OE-5}
\underline{T} V^* (p, i) \leq V^* (p, i), \qqd \forall p \in \pq(K), \ i \in S.
\de
Let $\mu_p$ be the measure given in (\ref{s3-OE-3}), $\llz_{\mu_p (\cdot | i), a} (p)$ is a probability measure on $K$ for any $i \in S$ and $a \in A$. By Theorem \ref{s3-thm3}, the value function $V^* (\llz_{\mu_p (\cdot | i), a} (p), i)$ exists, i.e.
\be\lb{s3-OE-6}
V^* (\llz_{\mu_p (\cdot | i), a} (p), i) = \overline{V} (\llz_{\mu_p (\cdot | i), a} (p), i) = \underline{V} (\llz_{\mu_p (\cdot | i), a} (p), i).
\de
The definition of $\underline{V} (\llz_{\mu_p (\cdot | i), a} (p), i)$ and (\ref{s3-OE-6}) imply that for each $\vz_0 > 0$, there exists some ${_{(i, a)}\pi} \in \ppz$, which depends on $i \in S$ and $a \in A$, such that
\be\lb{s3-OE-7}
V^* (\llz_{\mu_p (\cdot | i), a} (p), i) \leq V(\llz_{\mu_p (\cdot | i), a} (p), i, {_{(i, a)}\pz}, \sz) + \vz_0, \ \forall \sz \in \ssz.
\de
Define the policy $\pz^* = \{\pz_m^{*(k)}, k \in K, \ m \geq 0 \} \in \ppz$ as
\be\lb{s3-OE-add1}
\left\{
\begin{array}{ll}
\pz_0^{*(k)} (\cdot | i_0) =\mu_p^{(k)} (\cdot | i_0), & \hbox{$\forall i_0 \in S$;} \\
\pz_m^{*(k)} (\cdot | h_m) ={_{(i_0, a_0)}\pz}_{m-1}^{(k)} (\cdot | i_1,a_1,b_1, \ldots, i_{m}), & \hbox{$h_m=(i_0, a_0, b_0, \ldots, i_m)\in \Ho_m$.}
\end{array}
\right.
\de
Using the notation introduced in (\ref{1-pz}), one can verify that
\be\lb{s3-OE-8}
^{(i, a, b)}\pz_m^{*(k)} (\cdot | h_m) = \pz_{m+1}^{*(k)} (\cdot | i, a, b, h_m) = {_{(i,a)}\pz}_m^{(k)} (\cdot | h_m), \qd \text{$\forall k \in K$ and $(i,a,b) \in S\times A\times B$}.
\de
Hence, combining with (\ref{s3-lem7-0}), (\ref{s3-OE-7}) and (\ref{s3-OE-8}) , we obtain that
\begin{align*}
V & (p, i, \pz^*, \sz) \\
& = \sum_{k} \sum_{a, b} p_k \pz_0^{* (k)} (a | i) \sz_0 (b | i) c (k, i, a, b) \int_0^\infty e^{-\az t} (1 - D(t | i, a, b)) \d t \\
& \qd + \sum_{k} \sum_{a, b} p_k \pz_0^{* (k)} (a|i) \sz_0 (b|i) \sum_{i_1} \int_0^\infty e^{-\az t} Q (\d t, i_1| i, a, b) V (\llz_{\pz_0 (\cdot | i), a} (p), i_1, {_{(i, a)}\pz}, {^{(i, a, b)}\sz}) \\
& \geq \sum_{k } \sum_{a,b} p_k \mu_p ^{(k)} (a | i) \sz_0 (b | i) c (k, i, a, b) \int_0^\infty e^{-\az t} (1 - D(t | i, a, b)) \d t \\
& \qd +\sum_{k }  \sum_{a,b} p_k \mu_{p} ^{(k)} (a | i) \sz_0 (b|i) \sum_{i_1} \int_0^\infty e^{-\az t} Q (\d t, i_1 | i, a, b) \lt[ V^* (\llz_{\mu_{p} (\cdot | i), a} (p), i_1) - \vz_0 \rt]  \\
& \geq T^{\mu_{p}, \sz_0} V^* (p, i) - \vz_0.
\end{align*}
By the arbitrariness of $\sz \in \ssz$ and $\vz_0 >0$, (\ref{s3-OE-3}) implies that
$$
\underline{T} V^* (p, i) = \inf_{\sz_0 \in \ssz[0]} T^{\mz_{p}, \sz_0} V^* (p, i) \leq \underline{V} (p, i) = V^* (p, i).
$$

\underline{Step 3.} It can be proved in a similar way that
\be\lb{s3-OE-10}
V^* (p, i) \leq \overline{T} V^* (p, i), \qqd \forall p \in \pq(K), \ i \in \So.
\de
Here we only show the key steps. For arbitrary $\pz=\{\pz_n^{(k)}, k\in K, n\geq 0 \}\in \ppz$, denote $\pz_0=\{\pz_0^{(k)},k\in K\}\in \ppz[0]$, and then we obtain that $\llz_{\pz_0 (\cdot|i), a} (p)$ is a probability measure on $K$. Similar to (\ref{s3-OE-7}), for any $\vz_0>0$, $i\in S$ and $a\in A$, the definition of $\overline{V} (\llz_{\pz_{0} (\cdot | i), a} (p), i)$ ensures that there exists ${_{(i, a)}\sz} \in \ssz$ such that
\be\lb{s3-OE-11}
V^* (\llz_{\pz_{0} (\cdot | i), a} (p), i) \geq V(\llz_{\pz_0 (\cdot | i), a} (p),i , \pz, {_{(i, a)}\sz}) - \vz_0, \qd \forall \pz \in \ppz.
\de
Define the policy $\sz^* = \{\sz_m^*, m \geq 0 \} \in \ssz$ as
$$
\left\{
\begin{array}{ll}
\sz_0^{*} (\cdot | i) = \nu_p( \cdot | i ); & \\
\sz_m^{*} (\cdot | h_m) = {_{(i_0, a_0)}\sz}_{m-1}( \cdot |i_1, a_1, b_1, \ldots, i_m), & \hbox{$m \geq 1$, $h_m = ( i_0, a_0, b_0, \ldots, i_m)\in \Ho_m$,}
\end{array}
\right.
$$
where $\nu_p$ is given in (\ref{s3-OE-4}). It can be verified that
\be\lb{s3-OE-12}
^{(i, a, b)}\sz_m^{*} (\cdot | h_m) = \sz_{m+1}^{*} (\cdot |(i, a, b, h_m)) = {_{(i,a)}\sz}_m (\cdot | h_m), \qd \forall m \geq 0, \ h_m\in H_m.
\de
Again, combining with (\ref{s3-lem7-0}), (\ref{s3-OE-11}) and (\ref{s3-OE-12}), we obtain
\begin{align*}
V & (p, i, \pz, \sz^*) \\
& = \sum_{k} \sum_{a, b} p_k \pz_0^{(k)} (a | i) \sz^*_0 (b | i) c (k, i, a, b) \int_0^\infty e^{-\az t} (1 - D(t | i, a, b)) \d t \\
& \qd + \sum_{k} \sum_{a, b} p_k \pz_0^{(k)} (a|i) \sz^*_0 (b|i) \sum_{i_1} \int_0^\infty e^{-\az t} Q (\d t, i_1| i, a, b) V (\llz_{\pz_0 (\cdot | i), a} (p), i_1, {^{(i, a, b)}\pz}, {_{(i, a)}\sz}) \\
& \leq \sum_{k } \sum_{a, b } p_k \pz_0^{(k)} (a| i) \nu_{p} (b | i) c (k, i, a, b) \int_0^\infty e^{-\az t} (1 - D(t | i, a, b)) \d t \\
& \qd + \sum_{k } \sum_{a , b } p_k \pz_0^{(k)} (a | i) \nu_p(b|i) \sum_{i_1} \int_0^\infty e^{-\az t} Q (\d t, i_1 | i, a, b) \lt[ V^* (\llz_{\pz_0 (\cdot | i), a} (p), i_1) + \vz_0 \rt]  \\
\leq & T^{\pz_0, \nu_p} V^* (p, i) + \vz_0 .
\end{align*}
By the arbitrariness of $\pz \in \ppz$ and $\vz_0 >0$, we have
\be\lb{s3-OE-13}
V^* (p, i) = \underline{V} (p, i) \leq \sup_{\pz_0 \in \ppz[0]} T^{\pz_0, \nu_p} V^* (p, i) = \overline{T} V^* (p, i).
\de
Finally, (\ref{s3-OE-1}), (\ref{s3-OE-5}) and (\ref{s3-OE-13}) imply that $V^* = \overline{T} V^* = \underline{T} V^*$.
\deprf

In fact, the value function up to $n$-th decision epoch $V_n^* (p, i)$ satisfies the following iterative formula (\ref{s3-Vn-prop}). It should be noted that the order of limit and supremum (or infimum) in (\ref{VnOE}) cannot be exchanged, so Theorem \ref{OEthm} can not be proved by Corollary \ref{s3-Vn-prop}. On the contrary, the main proof method of Theorem \ref{OEthm} can be used to prove this corollary. Corollary \ref{s3-Vn-prop} also provides a method for iteratively calculating the value function $V_n^* (p, i)$.

\begin{cor}\lb{s3-Vn-prop}
For any $n\geq 0$, it holds that
\be\lb{VnOE}
V^*_{n+1}(p,i)=\overline{T}V^*_n(p,i)=\underline{T}V^*_n(p,i), \qd \forall p \in \pq(K), \ i \in S.
\de
Let $V_{-1}^* \equiv 0$, then we can calculate $V_n^*$ and $V^*$ by $V^*_{n} (p,i) = \sup_{\mu\in \ppz[0]} \inf_{\nu \in \ssz[0]} T^{\mu,\nu} V^*_{n-1} (p,i)$ and $V^*(p,i)=\lim_{n \to \infty}V^*_{n}(p,i)$.
\end{cor}

\prf The proof method is very similar to that of Theorem \ref{OEthm}, so we show the key steps and omit the details. Firstly, we can still use the Sion minimax theorem \cite[Theorem A.7]{SS02} to prove $\overline{T} V^*_n = \underline{T} V^*_n$. It is still valid to replace $V^*(p, i)$ with $V_n^*(p, i)$ in the proof method of Theorem \ref{OEthm}. Moreover, for each $p \in \pq (\Ko)$ there exist $\mu_{n, p} \in \ppz[0]$ and $\nu_{n, p} \in \ssz[0]$ such that
\begin{align}
\underline{T} V^*_n (p, i) &= \sup_{\mu \in \ppz[0]} \inf_{\nu \in \ssz[0]} T^{\mu, \nu} V^*_n (p, i) = \inf_{\nu \in \ssz[0]} T^{\mu_{n, p}, \nu} V^*_n (p, i), \qd \forall i\in S; \lb{s3-VnOE-3} \\
\overline{T} V^*_n (p, i) &= \inf_{\nu \in \ssz[0]} \sup_{\mu \in \ppz[0]} T^{\mu, \nu} V^*_n(p, i) = \sup_{\mu \in \ppz[0]} T^{\mu, \nu_{n, p}} V^*_n (p, i) \qd \forall i\in S. \lb{s3-VnOE-4}
\end{align}
Secondly, it is completely similar to the proof method of Step 2 of Theorem \ref{OEthm}. The key is the formula (\ref{s3-lem7-0-0}) in Lemma \ref{s3-lem7}, i.e.
\begin{align*}
V_{n+1} & (p, i, \pz^*, \sz) \\
& = \sum_{k} \sum_{a, b} p_k \pz_0^{* (k)} (a | i) \sz_0 (b | i) c (k, i, a, b) \int_0^\infty e^{-\az t} (1 - D(t | i, a, b)) \d t \\
& \qd + \sum_{k} \sum_{a, b} p_k \pz_0^{* (k)} (a|i) \sz_0 (b|i) \sum_{i_1} \int_0^\infty e^{-\az t} Q (\d t, i_1| i, a, b) V_n (\llz_{\pz_0 (\cdot | i), a} (p), i_1, {_{(i, a)}\pz}, {^{(i, a, b)}\sz}) \\
& \geq \sum_{k } \sum_{a,b} p_k \mu_{n, p}^{(k)} (a | i) \sz_0 (b | i) c (k, i, a, b) \int_0^\infty e^{-\az t} (1 - D(t | i, a, b)) \d t \\
& \qd +\sum_{k }  \sum_{a,b} p_k \mu_{n, p} ^{(k)} (a | i) \sz_0 (b|i) \sum_{i_1} \int_0^\infty e^{-\az t} Q (\d t, i_1 | i, a, b) \lt[ V_n^* (\llz_{\mu_{n, p} (\cdot | i), a} (p), i_1) - \vz_0 \rt]  \\
& \geq T^{\mu_{p}, \sz_0} V_n^* (p, i) - \vz_0,
\end{align*}
where $\mu_{n, p}$ is given in (\ref{s3-VnOE-3}) and $\pz^*$ is defined in a similar way as (\ref{s3-OE-add1}). Hence, we have
$$
\underline{T} V^*_n (p, i) = \inf_{\sz_0 \in \ssz[0]} T^{\mz_{n, p}, \sz_0} V^* _n(p, i) \leq V_{n+1}^* (p, i).
$$
Finally, we use the same method as Step 3 of Theorem \ref{OEthm} to prove that
$$
V^*_{n+1} (p, i) \leq \sup_{\pz_0 \in \ppz[0]} T^{\pz_0, \nu_{n, p}} V^*_n (p, i) = \overline{T} V^*_n (p, i).
$$
The above discussion can be concluded that $V^*_{n+1}(p,i)=\overline{T}V^*_n(p,i)=\underline{T}V^*_n(p,i)$.

\section{The existence of optimal policy for Player 1}
In this section, we focus on the existence of optimal policies for \pla 1 in $G(p)$ and its iterative algorithm. The main conclusion of this section is stated as follows.

\begin{thm}\lb{s4-thm1}
Suppose that Assumption \ref{ass1} holds. Given any $p\in \pq(\Ko)$ and $i \in S$, there exists $\pz^* \in \ppz$ such that
\be\lb{s4-thm1-1}
V(p, i, \pz^*, \sz) \geq V^*(p,i), \qd \forall \sz \in \ssz,
\de
i.e, $\pz^*$ is the optimal policy for \pla 1 in $G(p)$.
\end{thm}

Before proving Theorem \ref{s4-thm1}, we introduce a lemma to calculate the posterior probability of event $\{\kz = k\}$ given histories. This posterior probability also satisfies an iterative relation.

\begin{lem}\lb{s4-lem1}
For each $p \in \pq(K)$, $i \in \So$ and $(\pz, \sz) \in \ppz \times \ssz$. Denote by ${_{n}p}^{\pz}(k)$ the posterior distribution on $K$ at $n$-th decision epoch, i.e.,
$$
\left.
\begin{array}{ll}
	{_{0}p}^{\pz}(k) := \PP_{p,i}^{\pz,\sz} (\kz=k | H_0 = i_0, T_{0} = 0) = p_k, & \\
	{_{n}p}^{\pz}(k) := \PP_{p,i}^{\pz,\sz} (\kz=k | H_m = h_m, \ddz T_{m} \leq s_m, 1 \leq m \leq n), & \hbox{$ \forall n \geq 1$,}
\end{array}
\right.
$$
where $\ddz T_m := T_m -T_{m-1}$. Then, ${_{n}p}^{\pz}(k)$ satisfies that
$$
{_{n+1}p}^{\pz} (k) = \llz_{\pz_n(\cdot | h_n), a_n} ({_{n}p}^\pz) (k), \qd \forall k \in K, \ n \geq 0.
$$
\end{lem}

\prf Using the expressions (\ref{sect2-1}), (\ref{sect2-2}) and the Bayes formula, we obtain
\begin{align*}
{_{n}p}^{\pz}(k) & = \PP_{p,i}^{\pz,\sz} (\kz=k | H_m=h_m, \ddz T_{m} \leq s_m, 1 \leq m \leq n) \\
& = \frac{p_k \lt[\prod_{m=0}^{n-1}\pz_{m}^{(k)}(a_m|h_m)\sz_{m}(b_m|h_m) \int_0^{s_n} Q(\d t, i_m| i_{m-1}, a_{m-1}, b_{m-1}) \rt]}{\sum_{l \in K} p_l \lt[\prod_{m=0}^{n-1} \pz_{m}^{(l)} (a_m| h_m) \sz_{m} (b_m| h_m) \int_0^{s_m} Q(\d t, i_m| i_{m-1}, a_{m-1}, b_{m-1}) \rt]} \\
& = \frac{p_k \prod_{m=0}^{n-1} \pz_m^{(k)} (a_m| h_m)}{\sum_{l\in K} p_l \prod_{m=0}^{n-1} \pz_m^{(l)} (a_m| h_m)}.
\end{align*}
Obviously, it holds that ${_{n}p}^{\pz} \in \pq(K)$. Moreover, ${_{n+1}p}^{\pz}$ can be represented by ${_{n}p}^{\pz}$, and then using the definition of $\llz$ introduced in (\ref{llz}), we have
$$
{_{n+1}p}^\pz (k) = \frac{ {_{n}p}^\pz(k) \pz_{n}^{(k)}(a_{n}|h_{n})}{\sum_{l\in K}{_{n}p}^\pz(l) \pz_{n}^{(l)} (a_{n} | h_{n})} = \llz_{\pz_{n}(\cdot|h_{n}),a_{n}} ({_{n}p}^\pz) (k), \qd \forall k \in K, \ n \geq 0,
$$
where $\pi_{n} ( \cdot | h_{n}) = \{\pz_{n}^{(k)} (\cdot | h_{n}), k \in K \} \in \pq(A|K)$. This leads to the required recurrence relation.
\deprf

This is consistent with the intuition, that is, the posterior distribution on $K$ at $n$-th decision epoch has nothing to do with the actions of \pla 2, but only depends on \pla 1. The essential reason is that \pla 2 has asymmetric information and do not know the selection of the system.

{\bf Proof of Theorem \ref{s4-thm1}.} Based on Lemma \ref{s4-lem1}, we can construct the optimal policy for \pla 1 step by step. Firstly, we define the pairs of measures $\lt\{({_{n}p}^*, \pz_n^{*} (\cdot | h_n)), n \geq 0 \rt\}$ recursively. That is
\be\lb{s4-thm1-4}
\left\{
\begin{array}{ll}
	{_{0}p}^* (k) := p_k, & \hbox{$k \in K$;} \\
	{_{n+1}p}^* (k) := \llz_{\mz[_{n}p^*] (\cdot | i_n), a_n} (_{n}p^*) (k), & \hbox{$n \geq 0$ and $k \in K$,}
\end{array}
\right.
\de
where $\mz[_{n}p^*]$ is constructed by $_{n}p^*$ according to (\ref{s3-OE-3}) in Step 1 of the proof of Theorem \ref{OEthm}. In detail, by induction, we have $_{n}p^* \in \pq(K)$ for all $n \geq 0$. Then, for each $_{n}p^*$, using Theorem \ref{OEthm} and (\ref{s3-OE-3}), there is $\mz[_{n}p^*] \in \ppz[0]$ such that
\be\lb{s4-thm1-5}
V^* (_{n}p^*, i) = \sup_{\mz \in \ppz[0]} \inf_{\nz \in \ssz[0]} T^{\mz, \nz} V^* (_{n}p^*, i) = \inf_{\nz \in \ssz[0]} T^{\mz[_{n}p^*], \nz} V^* (_{n}p^*, i), \qd \forall i\in S.
\de
Denote by
\be\lb{s4-thm1-6}
\pz_n^{*(k)} (\cdot | h_n) = \mz[_{n}p^*]^{(k)} (\cdot | i_n), \qqd \forall k \in K.
\de
Hence, for each $n \geq 0$, we have ${_{n}p}^* \in \pq(K)$ and $\pi^* := \{\pi^{*(k)}_n, k\in K, n\geq 0\} \in \Pi$. Furthermore, Lemma \ref{s4-lem1} implies that
\be\lb{s4-thm1-7}
\PP_{p,i}^{\pz,\sz} (\kz=k | H_m = h_m, \Delta T_{m} \leq s_m, 0 \leq m \leq n+1) = {_{n+1}p}^{*} = \llz_{ \mz[_{n}p^*] (\cdot | i_n), a_n} ({_{n}p}^{*}).
\de

Next, for every fixed $\sigma=\{\sigma_n,n\geq 0\}\in \Sigma$, we claim that for each $n \geq 0$,
\begin{align}\lb{s4-thm1-8}
V^* (p, i) & \leq \frac{1}{\az} \sum_{m=0}^{n} \EE_{p, i}^{\pz^*, \sz} \lt[(e^{-\az T_m} - e^{-\az T_{m+1}}) c(\kz, X_m, A_m, B_m) \rt] \notag \\
& \qd + \EE_{p, i}^{\pz^*, \sz} \lt[ e^{-\az T_{n+1}} V^* ({_{n+1}p}^*, X_{n+1}) \rt].
\end{align}
This inequality can be proved by induction. Clearly, it holds that for $n=0$, since
\begin{align*}
V^*(p,i) & \leq T^{\pi^*_0,\sigma_0}V^*(p,i)\\
& = \frac{1}{\alpha}\EE^{\pi^*,\sigma}_{p,i}\lt [(1-e^{\alpha T_1})c(\kz, X_0, A_0, B_0)\rt]\\
& \qd +\sum_{k}\sum_{a, b} p_k \pi^{*(k)}_0(a|i)\sigma_0(b|i)\sum_{j\in S}\int_0^\infty Q(\d t, j | i, a, b) e^{-\alpha t} V^*({_1p}^{*}, j) \\
& = \frac{1}{\alpha} \EE^{\pi^*, \sigma}_{p,i} \lt[ (1 - e^{\alpha T_1}) c(\kz, X_0, A_0, B_0)\rt] + \EE^{\pi^*, \sigma}_{p,i} \lt[ e^{-\alpha T_1} V^*({_1p}^*, X_1) \rt],
\end{align*}
where the first inequality and the seconde equality hold according to the equations (\ref{s4-thm1-5}) and (\ref{s4-thm1-7}) respectively.

Next, we consider $n \geq 1$. Given any $h_n \in \Ho_n$, $a_n \in A$ and $b_n \in B$, we have $\sigma_{n+1} (\cdot | h_{n}, a_{n}, b_{n}, \cdot) \in \Sigma[0]$, and then
\begin{align*}
& V^*({_{n+1}p}^*, i_{n+1}) \leq T^{ \mu[ {_{n+1}}p^{*} ] , \sigma_{n+1} (\cdot | h_{n}, a_{n}, b_{n}, \cdot)} V^*({_{n+1}}p^{*}, i_{n+1})\\
& = \Bigg\{\sum_{k} {_{n+1}p}^* (k) \sum_{a, b} \pi^{*(k)}_{n+1} (a | h_{n+1}) \sigma_{n+1} (b | h_{n+1}) c(k, i_{n+1}, a, b) \int_0^\infty e^{-\alpha t} (1-D( t|i_{n+1},a,b))\d t \Bigg\} \\
& \qd + \Bigg\{ \sum_{k} {_{n+1}p}^* (k) \sum_{a, b} \pi^{*(k)}_{n+1} (a | h_{n+1}) \sigma_{n+1} ( b | h_{n+1})  \sum_{j} \int_0^\infty e^{-\alpha t} Q(\d t,j | i_{n+1}, a, b) V^* ({_{n+2}p}^{\pi^*}, j) \Bigg\} \\
& = \frac{1}{\alpha} \EE_{p,i}^{\pi^*,\sigma} \lt[\lt(1 - e^{-\alpha (T_{n+2} - T_{n+1})}\rt) c(\kz, X_{n+1}, A_{n+1}, B_{n+1}) \big| H_m = h_m, \Delta T_{m} \leq s_m, 0 \leq m \leq n+1 \rt] \\
& \qd + \EE_{p,i}^{\pi^*, \sigma} \lt[ e^{-\alpha (T_{n+2} - T_{n+1})} V^*({_{n+2}p^*(H_{n+1}, A_{n+1})}, X_{n+2}) \big| H_m = h_m, \Delta T_{m} \leq s_m, 0 \leq m \leq n+1 \rt].
\end{align*}
Hence, using (\ref{sect2-1}) and (\ref{sect2-2}), the second item of (\ref{s4-thm1-8}) satisfies
\begin{align}\lb{s4-thm1-9}
& \EE_{p,i}^{\pi^*,\sigma} \lt[ e^{-\alpha T_{n+1}}V^*({_{n+1}p}^*, X_{n+1}) \rt] \notag \\
& \leq \frac{1}{\alpha} \EE_{p,i}^{\pi^*,\sigma} \lt[ \lt( e^{-\alpha T_{n+1}} - e^{-\alpha T_{n+2}} \rt) c(\kz, X_{n+1}, A_{n+1}, B_{n+1}) \rt] + \EE_{p,i}^{\pi^*,\sigma} \lt[ e^{-\alpha T_{n+2}} V^*({_{n+2}p^*}, X_{n+2}) \rt].
\end{align}
Assume that (\ref{s4-thm1-8}) holds for $N=n$. For $N= n+1$, by (\ref{s4-thm1-9}), we have
\begin{align*}
V^* (p, i) & \leq \frac{1}{\az} \sum_{m=0}^{n} \EE_{p, i}^{\pz^*, \sz} \lt[ \lt( e^{-\az T_m} - e^{-\az T_{m+1}} \rt) c(\kz, X_m, A_m, B_m) \rt] + \EE_{p, i}^{\pz^*, \sz} \lt[ e^{-\az T_{n+1}} V^* ({_{n+1}p}^*, X_{n+1}) \rt] \\
& \leq \frac{1}{\az} \sum_{m=0}^{n+1} \EE_{p, i}^{\pz^*, \sz} \lt[ \lt( e^{-\az T_m} - e^{-\az T_{m+1}} \rt) c(\kz, X_m, A_m, B_m) \rt] + \EE_{p,i}^{\pi^*,\sigma}\lt [e^{-\alpha T_{n+2}}V^*({_{n+2}p^*}, X_{n+2}) \rt].
\end{align*}
This completes the induction proof of the inequality (\ref{s4-thm1-8}). Using (\ref{s2-prop2-2}), we have $V^* (p, i) \leq V (p, i, \pz^*, \sz)$ by passing $n \to \infty$ in (\ref{s4-thm1-8}). Finally, by the arbitrariness of $\sz \in \ssz$ and $i \in S$, we obtain the existence of the optimal policy for \pla 1 in $G(p)$.
\deprf

The proof of Theorem \ref{s4-thm1} also gives the iterative algorithm of the optimal policy for \pla 1. The key point is that the probability distribution of game type and the optimal policy for \pla 1 need to be calculated together, i.e., $\lt\{({_{n}p}^*, \pz_n^{*} (\cdot | h_n)), n \geq 0 \rt\}$. When iteratively calculating the optimal policy, the original game $G(p)$ is regarded as a new game $G(_{n}{p}^*)$ due to the change of probability distribution $_{n}{p}^*$ at the $n$-th decision epoch. For the convenience of application, the methods in the proof are arranged into the following algorithms.

\begin{breakablealgorithm}\lb{algo-p1}
\caption{optimal policy for \pla 1} 
\hspace*{0.2in} {\bf Input:} 
The two-players zero-sum semi-Markov game with incomplete information $\{\Ko, \So, (\Ao \times \Bo), p, Q(\cdot, \cdot|i, a, b), c(k, i, a, b) \}$; the value function $V^*$ given by Theorem \ref{s3-thm3} and Corollary \ref{s3-Vn-prop}; for each $n \geq 0$, the history $h_n = (i_0, a_0, b_0, \ldots, i_n) \in \Ho_n$.
\begin{algorithmic}[1]
\State Let ${_{0}p^*} (k) := p_k$, $k \in K$.
\For{$n=0, 1, 2, \ldots $} \\ 
Compute
$$
\mu[ {_{n}p^*} ] := \arg \max_{\mu \in \ppz[0]} \{\min_{\nu \in \ssz[0]} T^{\mu,\nu} V^*({_{n}p^*}, i_n) \}.
$$
Update the police at the $n$-th decision epoch $\pz_{n}^{*(k)} (\cdot|h_{n}) := \mu[{_{n}p^*}]$. \\
Compute
$$
{_{n+1}p}^{*} (k) = \llz_{ \mz[_{n}p^*] (\cdot | i_n), a_n} ({_{n}p}^{*}) (k) = \frac{p_k \prod_{m=0}^n \pz_m^{*(k)} (a_m|h_m)}{\sum_{l\in K} p_l \prod_{m=0}^n \pz_m^{*(l)} (a_m|h_m)}.
$$
Update the probability distribution of game type at the $(n+1)$-th decision epoch by ${_{n+1}p}^{*}$ and repeat step 3.
\EndFor
\end{algorithmic}
\hspace*{0.2in} {\bf Output:} 
The policy $\pz^*:=\{\pz^{*(k)}_n, k\in K, n\geq 0\}$ is the optimal policy for \pla 1.
\end{breakablealgorithm}

\section{The dual games and the existence of the optimal policy for Player 2}

In this section, we focus on the optimal policy for \pla 2. To do so, we introduce the concept of dual games with incomplete information. The existence of the optimal policy for \pla 2 in $G(p)$ can be proved by the relationship between the dual game and the original game. The dual semi-Markov game with incomplete information is the six-tuple given as:
$$
\{K, S, A \times B, z, Q(\cdot, \cdot |i, a, b), c(k, i, a, b) \},
$$
where $K$, $S$, $A\times B$, $Q$ and $c$ are the same as the original game $G(p)$ defined in Section 2. The difference is that the probability $p \in \pq(K)$ in $G(p)$  is replaced by the real-value vector $z = \{z(k), k \in K\} \in \RR^{|K|}$, which is used to modify the cost in the value function (see (\ref{dual_val}) below). In the dual game, the type of game $\kz = k$ is not determined by the system, but is choosed by \pla 1. Instead, the cost function is modified by the real-value function $z(k)$. Denote by $G^{\#} (z)$ the dual semi-Markov game with incomplete information on one side.

In detail, the dual game $G^{\#} (z)$ evolves in the following way. At the initial decision time $t_0 = 0$, the system stays at $i_0$ and \pla 1 chooses a game type $k \in K$ according to the initial state $i_0$. The game type will not change in the subsequent evolution, but is hidden from \pla 2. In the meantime, \pla 1 and \pla 2 choose action $a_0 \in A$ and $b_0 \in B$ respectively base on the current state $i_0$. The actions chosen are perfectly observed by both players. Similar to the original game $G(p)$, the system stays at $i_0$ no more than time $t_1$ and jumps to $i_1$ according to the semi-Markov kernel $Q (t_1, i_1 | i_0, a_0, b_0)$. Then, the next decision epoch occurs. At this moment \pla 1 chooses an action $a_1 \in A$ base on game type $k$ and $(i_0, a_0, b_0, i_1)$, \pla 2 chooses an action $b_1 \in B$ base on $(i_0, a_0, b_0, i_1)$. The system evolves repeatedly in the above way.

The policy for \pla 2 in the dual game is coincide with that in the original game, refer to Definition \ref{pol} in Section 2. Nevertheless, the policy for \pla 1 in the dual game $G^{\#} (z)$ is different, which is given below.

\begin{defn}
A randomized history-dependent policy $\hat{\pz}$ for \pla 1 in the dual game $G^{\#} (z)$ is given by two-tuples $(\rz, \pz)$, where $\rz = \{ \rz_{i} ,i \in S\} \in \pq(K|S)$ and $\pz = \{ \pz_n^{(k)}, k \in K, n \geq 0\} \in \ppz$. Denote by $\hat{\ppz}$ the set of all policies for \pla 1 in the dual game $G^{\#} (z)$, i.e., $\hat{\ppz}=\pq(K|S)\times \ppz$.
\end{defn}

In the dual game $G^{\#} (z)$, the notations of the system are coincide with that in the original game $G(p)$, which include the trajectory space $(\ooz, \fq)$, history $H_n$ up to the $n$-th decision epoch and the random variables $\kz$, $X_n$, $T_n$, $A_n$ and $B_n$. For each $i \in S$, $\hat{\pz} =(\rz,\pz) \in \hat{\ppz}$ and $\sz \in \ssz$, Tulcea's Theorem (\cite[Proposition C.10]{HL96}) implies that there exist a unique probability measure $\PP_{i}^{\hat{\pz}, \sz}$ on $(\ooz, \fq)$ satisfying
\begin{align}
&\PP_{i}^{\hat{\pz},\sz}(X_0=i,\kz=k)=\rz_i(k) \lb{dual-p1}\\
&\PP_{i}^{{\hat{\pz},\sz}}(A_n = a, B_n=b | \kz, H_n, T_n)=\pz_n^{(\kz)} (a | H_n) \sz_{n} (b | H_n) \lb{dual-p2}\\
&\PP_{i}^{{\hat{\pz},\sz}}(T_{n+1}-T_{n}\leq t,X_{n+1}=j | \kz,H_n, T_n, A_n, B_n) = Q ( t, j | X_n, A_n, B_n) \lb{dual-p3}
\end{align}
According to the relationship between the dual game and the original game, we have
\be\lb{dual_eq1}
\PP_{i}^{\hat{\pz}, \sz} (E) = \PP_{\rz_i, i}^{\pz, \sz} (E), \qd \forall E \in \fq.
\de
Denote by $\EE_{i}^{\hat{\pz}, \sz}$ the expectation corresponding to $\PP_{i}^{\hat{\pz}, \sz}$. In the dual game $G^{\#} (z)$, for each $\hat{\pz} =(\rz,\pz)\in \hat{\ppz}$ and $\sz \in \ssz$, the expected discount reward for \pla 1 is defined as
\begin{align}\lb{dual_val}
U (z, i, \hat{\pz}, \sz)
:= &\EE_i^{\hat{\pz}, \sz} \lt[ \int_0^{T_{\infty}} e^{-\az t} \lt( c(\kz, X_t, A_t, B_t) -  z(\kz) \rt) \d t \rt]\notag\\
= &\sum_{m=0}^\infty \EE_i^{\hat{\pz}, \sz} \lt[ \frac{1}{\az} \lt(e^{-\az T_m} - e^{-\az T_{m+1}} \rt) c(\kz, X_m, A_m, B_m) \rt] -  \EE_i^{\hat{\pz}, \sz} \lt [z(\kz) \int_0^{T_\infty} e^{-\az t}\d t\rt].
\end{align}
Under Assumption \ref{ass1}, by (\ref{dual_eq1}), we have $\PP_{i}^{\hat{\pz}, \sz} \lt( T_{\infty} = \infty \rt) =1$ and
\begin{align}\lb{dual_eq2}
U (z, i, \hat{\pz}, \sz) &= \EE_{\rz_i, i}^{\pz, \sz} \lt[ \int_0^\infty e^{-\az t} c (\kz, X_t, A_t, B_t) \d t \rt] - \frac{1}{\az} \sum_{k \in K} z(k) \PP_{\rz_i, i}^{\pz, \sz} (\kz = k) \notag \\
&= V(\rz_i, i, \pz, \sz) - \frac{1}{\az} \lan \rz_i, z \ran,
\end{align}
where $\lan z, \rz_i \ran := \sum_{k \in K} z(k) \rz_i (k)$.
Similar to the original game $G(p)$, we give the definitions of lower value $\underline{U} (z, i)$ and upper value $\overline{U} (z, i)$ as
\begin{align}
\underline{U} (z, i) &= \sup_{\hat{\pz} \in \hat{\ppz}} \inf_{\sz \in \ssz} U (z, i, \hat{\pz}, \sz)=\sup_{\rho\in \pq(K|S)}\sup_{\pz\in \ppz} \inf_{\sz \in \ssz} U (z, i, (\rho,\pz), \sz) \lb{d-s}\\
\overline{U} (z, i) &= \inf_{\sz \in \ssz} \sup_{\hat{\pz} \in \hat{\ppz}} U (z, i, \hat{\pz}, \sz) = \inf_{\sz \in \ssz}\sup_{\rho\in \pq(K|S)}\sup_{\pz\in \ppz} U (z, i, (\rho,\pz), \sz) \lb{d-i}.
\end{align}
If $\underline{U} = \overline{U}$ holds, we say the value function of the dual game $G^{\#} (z)$ exists and denote
$$U^* (z, i) = \underline{U} (z, i) = \overline{U} (z, i).$$
Here, we give the definition of the optimal policy in the dual game $G^{\#} (z)$, which is similar to Definition \ref{optimal}.
\begin{defn}\lb{optimal0}
Fix any $z \in \RR^{|K|}$. A policy $\hat{\pz}^* \in \hat{\ppz}$ for \pla 1 is called optimal in the dual game $G^{\#} (z)$ if
$$
\inf_{\sz \in \ssz} U (z, i, \hat{\pz}^*, \sz) \geq \underline{U} (z, i), \qd \forall i\in S.
$$
A policy $\sz^* \in \ssz$ for \pla 2 is called optimal in the dual game $G^{\#} (z)$ if
$$
\sup_{\hat{\pz} \in \hat{\ppz}} U(z, i, \hat{\pz}, \sz^*) \leq \overline{U} (z, i)\qd \forall i\in S.
$$
\end{defn}

Similar to the value function of the original game $G(p)$, $\underline{U}$ or $\overline{U}$ has continuity (Lemma \ref{s5-lem1}) and convexity (Lemma \ref{s5-lem3}). These properties are described later here.

\begin{lem}\lb{s5-lem1}
For each $i \in S$, $\underline{U} (z, i)$ and $\overline{U} (z, i)$ are Lipschitz continuous functions respect to $z \in \RR^{|K|}$.
\end{lem}
\prf Given any $(\hat{\pz} ,\sz)\in \hat{\ppz} \times \ssz$ and $z_1=\{z_{1} (k),k\in K\}, z_{2} = \{z_{2} (k), k\in K\} \in \RR^{|K|}$, we have
$$
\lt| U (z_1, i, \hat{\pz}, \sz) - U (z_2, i, \hat{\pz}, \sz) \rt| \leq \EE_i^{\hat{\pz}, \sz} \lt[ \frac{1}{\az} \lt| z_1 (\kz) - z_2 (\kz) \rt| \rt] \leq \frac{1}{\az} \| z_1 - z_2 \|,
$$
The rest is similar to the proof of Lemma \ref{s3-lem4}.
\deprf

\begin{lem}\lb{s5-lem3}
	Suppose that Assumption \ref{ass1} holds. For each $i \in S$, $\overline{U} (\cdot, i)$ is a convex function on $\RR^{|K|}$.
\end{lem}

\prf Fix arbitrary $i \in S$ and $\vz > 0$. For any $z_1, z_2 \in \RR^{|K|}$, by the definition of $\overline{U}$, there exists $\sz_1, \sz_2 \in \ssz$ such that
$$
\sup_{\hat{\pz} \in \hat{\ppz}} U (z_m, i, \hat{\pz}, \sz_m) \leq \overline{U} (z_m, i) + \vz, \qd m=1, 2.
$$
Using Proposition \ref{s3-thm5} and the equation (\ref{dual_eq1}), for each $\lz \in [0, 1]$, there exists $\sz^\lz \in \ssz$ such that
\begin{align*}
V(p,i,\pz,\sz^\lz)=\lz V(p,i,\pz,\sz_1)+(1-\lz)V(p,i,\pz,\sz_2),\qd \forall i\in S,p\in \pq(K),\pz\in \ppz.
\end{align*}
For each $\hat{\pz} = (\rho, \pz) \in \hat{\ppz}$, using Assumption \ref{ass1}, we have
\begin{align*}
 &\lz U(z_1, i, \hat{\pz}, \sz_1) + (1- \lz) U(z_2, i, \hat{\pz}, \sz_2) \\
& = \lz V(\rho_i,i,\pz,\sz_1)+(1-\lz)V(\rho_i,i,\pz,\sz_2) - \frac{\lz}{\az}\lan \rho_i,z_1\ran -\frac{1-\lz}{\az} \lan \rho_i,z_2\ran\\
&= V(\rho_i,i,\pz,\sz^\lz) - \frac{1}{\az} \lan \lz z_1 + (1 - \lz) z_2, \rz_i \ran \\
&= U(\lz z_1 + (1 - \lz) z_2, i, \hat{\pz},  \sz^\lz).
\end{align*}
By the arbitrariness of $\hat{\pz}$, we have
\begin{align*}
\overline{U} (\lz z_1 + (1- \lz) z_2, i) & \leq \sup_{\hat{\pz} \in \hat{\ppz}} U (\lz z_1 + (1- \lz) z_2, i, \hat{\pz},  \sz^\lz) \\
& \leq \lz \sup_{\hat{\pz} \in \hat{\ppz}} U (z_1, i, \hat{\pz}, \sz_1) + (1 - \lz) \sup_{\hat{\pz} \in \hat{\ppz}} U (z_2, i, \hat{\pz}, \sz_2) \\
& \leq \lz \overline{U} (z_1, i) + (1- \lz) \overline{U} (z_2, i) + \vz,
\end{align*}
which, together with the arbitrariness of $\vz$,  means that $\overline{U} (\cdot, i)$ is convex on $\RR^{|K|}$.\deprf

The next results are about the relationship between the value function of the original game $G(p)$ and the one of the dual game $G^{\#} (z)$, which are given by the variational expressions.

\begin{thm}\lb{s5-thm4}
	Suppose that Assumption \ref{ass1} holds. For each $i \in S$, we have
	\begin{align}
	\underline{U} (z, i) = \max_{p \in \pq (K)} \lt\{ \underline{V} (p, i) - \frac{1}{\az} \lan p, z \ran \rt\}, \qd z \in \RR^{|K|}; \lb{s5-thm4-1} \\
	\overline{U} (z, i) = \max_{p \in \pq (K)} \lt\{ \overline{V} (p, i) - \frac{1}{\az} \lan p, z \ran \rt\}, \qd z \in \RR^{|K|}. \lb{s5-thm4-2}
	\end{align}
	In the dual case, the lower and upper value functions of the game $G(p)$ satisfy
	\begin{align}
	\underline{V} (p, i) = \min_{z \in \mathbb{B}} \lt\{ \underline{U} (z, i) + \frac{1}{\az} \lan p, z \ran \rt\}, \qd p \in \pq (K); \lb{s5-thm4-3} \\
	\overline{V} (p, i) = \min_{z \in \mathbb{B}} \lt\{ \overline{U} (p, i) + \frac{1}{\az} \lan p, z \ran \rt\}, \qd p \in \pq (K), \lb{s5-thm4-4}
	\end{align}
	where $\mathbb{B} := \lt\{ z\in\RR^{|K|} \Big| 0\leq z(k) \leq c^*, \forall k \in K \rt\}$ is a compact subset of $\RR^{|K|}$.
\end{thm}

\prf In the following, we give the direct proof of (\ref{s5-thm4-1}) and (\ref{s5-thm4-4}), while for the other two formulas (\ref{s5-thm4-2}) and (\ref{s5-thm4-3}) are proved by the Fenchel theorem given in \cite[Theorem A. 16]{SS02}. The following proof is divided into three steps.

\underline{Step 1.} According to (\ref{dual_eq2}), for each $z \in \RR^{|K|}$ and $i \in S$, we have
\begin{align*}
\underline{U} (z, i) &= \sup_{\rho \in \pq(K|S)} \sup_{\pz\in \ppz} \inf_{\sz \in \ssz} \lt\{  V(\rho_i,i,\pz,\sz)- \frac{1}{\az} \lan \rz_i, z \ran \rt\} \\
&= \sup_{\rz_i \in \pq (K)}\lt\{ \sup_{\pz \in \ppz} \inf_{\sz \in \ssz}   V(\rho_i,i,\pz,\sz)- \frac{1}{\az} \lan \rz_i, z \ran \rt\} \\
&= \sup_{p \in \pq (K)} \lt\{ \underline{V} (p, i) - \frac{1}{\az} \lan p, z \ran \rt\}.
\end{align*}
By the Lemma \ref{s3-lem4} and compactness of $\pq (K)$, the supremum can be replaced by the maximum in above equation, i.e.,
$$
\underline{U} (z, i)=  \max_{p \in \pq (K)} \lt\{ \underline{V} (p, i) - \frac{1}{\az} \lan p, z \ran \rt\}.
$$

\underline{Step 2.} In this step, we give the proof of (\ref{s5-thm4-4}). Fix any $i \in S$ and $z \in \RR^{|K|}$. Note that for any $\vz > 0$, there exists a policy $\widetilde{\sz} = \widetilde{\sz} (z, i, \vz)$ such that
$$
\sup_{\pz \in \ppz} \EE_{p, i}^{\pz, \widetilde{\sz}} \lt[ \int_0^\infty e^{-\az t} \lt( c(\kz, X_t, A_t, B_t) - z(\kz) \rt) \d t \rt] \leq \overline{U} (z, i) + \vz, \qd \forall p \in \pq (K).
$$
Hence, the definition of $\overline{V} (p, i)$ implies that
$$
\overline{V} (p, i) \leq \sup_{\pz \in \ppz} V(p,i,\pi,\widetilde{\sz}) \leq  \overline{U} (z, i) + \frac{1}{\az} \lan p, z \ran  + \vz, \qd \forall p \in \pq (K).
$$
By the arbitrariness of $z$ and $\vz$, we have that
\be\lb{s5-thm4-6}
\overline{V} (p, i)  \leq  \inf_{z\in \RR^{|K|}}\lt\{  \overline{U} (z, i) + \frac{1}{\az} \lan p, z \ran \rt \} \leq  \inf_{z\in \mathbb{B}}\lt\{  \overline{U} (z, i) + \frac{1}{\az} \lan p, z \ran \rt \}, \qd \forall p \in \pq (K).
\de
Conversely, fix arbitrary $i \in S$ and $p \in \pq (K)$. By the definition of $\overline{V} (p, i)$, for any $\vz >0$, there exists $\hat{\sz} = \hat{\sz} (p, i, \vz)$ such that
\begin{align}\lb{s5-thm4-7-1}
\overline{V} (p, i) + \vz & \geq \sup_{\pz \in \ppz}V(p,i,\pz,\hat{\sz})  = \sum_{k \in K} p_k \sup_{\pz \in \ppz} V(\dz_k,i,\pz,\hat{\sz}).
\end{align}
We need to show that second equality of (\ref{s5-thm4-7-1}) holds. Fix $p = \dz_k$, $k \in K$. According to (\ref{E[f]}), the value of $\EE_{\dz_k, i}^{\pz, \hat{\sz}} [f]$ only depends on $\{ \pz_n^{(k)} \}_{n \geq 0}$, i.e., it is independent of $\{ \pz_n^{(j)} \}_{n \geq 0}$, $(j \neq k)$. Then, using Theorem \ref{s4-thm1}, the optimal policy of \pla 1 exists and is denoted by $\pz^{* (k)} = \{ \pz_n^{* (k)}, n \geq 0 \}$. Let $\pz^* := \{\pz^{* (k)}, k \in K \}$ and using (\ref{E[f]}) again, we have
\begin{align*}
\sum_{k \in K} p_k \sup_{\pz \in \ppz}  V(\dz_k,i,\pz,\hat{\sz})& = \sum_{k\in K} p_k V(\dz_k,i,\pz^*,\hat{\sz})  \leq\sup_{\pz \in \ppz} V(p,i,\pz,\hat{\sz}).
\end{align*}
The inverse inequality is obvious, and then the second equality of (\ref{s5-thm4-7-1}) holds. Define $z^{\hat{\sz}}=\{z^{\hat{\sz}} (k), k\in K\} $ as
\be\lb{z}
z^{\hat{\sz}} (k) := \az \sup_{\pz \in \ppz} \EE_{\dz_k, i}^{\pz, \hat{\sz}} \lt[ \int_0^\infty e^{-\az t} c (\kz, X_t, A_t, B_t) \d t \rt]=
\az  \sup_{\pz \in \ppz}  V(\dz_k,i, \pz,\hat{\sz}) .
\de
Obviously, it holds that $0 \leq z^{\hat{\sz}} (k) \leq c^*$ for every $k\in K$, i.e., $z^{\hat{\sz}}\in \mathbb{B}$. In the dual game $G(z^{\hat{\sz}})$, we have
\begin{align*}
\overline{U} (z^{\hat{\sz}}, i) & \leq \sup_{p \in \pq(K)} \sup_{\pz \in \ppz} \lt\{V(p,i,\pz,\hat{\sz}) - \frac{1}{\az} \lan p, z^{\hat{\sz}} \ran \rt\}\\
& = \sup_{p \in \pq (K)} \sup_{\pz \in \ppz}\lt\{ \sum_{k\in K}p_k \lt( V(\dz_k,i,\pz,\hat{\sz}) - \frac{1}{\az} z^{\hat{\sz}} (k) \rt) \rt\} \leq 0,
\end{align*}
where the last inequality is based on (\ref{z}). Then, combining (\ref{s5-thm4-7-1}) and (\ref{z}), we obtain
\begin{align}\lb{s5-thm4-7}
\overline{V} (p, i) & \geq \frac{1}{\az} \lan p, z^{\hat{\sz}} \ran + \overline{U} (z^{\hat{\sz}}, i) - \vz \notag \\
& \geq \inf_{z \in \mathbb{B}} \lt\{ \frac{1}{\az} \lan p, z \ran + \overline{U} (z, i) \rt\} - \vz
\geq \inf_{z \in \RR^{|K|}} \lt\{ \frac{1}{\az} \lan p, z \ran + \overline{U} (z, i) \rt\} - \vz.
\end{align}
According to (\ref{s5-thm4-6}), (\ref{s5-thm4-7}) and the arbitrariness of $\vz$, we obtain
\be\lb{s5-thm4-c-1}
\overline{V} (p, i) = \inf_{z \in \RR^{|K|}} \lt\{ \frac{1}{\az} \lan p, z \ran + \overline{U} (z, i) \rt\}=\inf_{z \in \mathbb{B}} \lt\{ \frac{1}{\az} \lan p, z \ran + \overline{U} (z, i) \rt\}.
\de
According to the compactness of $\mathbb{B}$ and the continuity of $\overline{U}$ given in Lemma \ref{s5-lem1}, the infimum can be replaced by the minimum, which is (\ref{s5-thm4-4}).

\underline{Step 3.} In order to show (\ref{s5-thm4-2}) and (\ref{s5-thm4-3}), we introduce the concepts of the Fenchel duality. Let $f$ be a function defined on $\RR^n$ with values in $\RR \cup \{+ \infty \}$. The Fenchel conjugate of $f$ is the function defined on $\RR^n$
\be\lb{s5-thm4-5}
f^{\#} (p) := \sup_{x \in \RR^n} \{ \lan x, p \ran - f(x) \},
\de
where $\lan \cdot, \cdot \ran$ is the inner product on $\RR^n$. Note that $\pq (K)$ is a subset of $\RR^{|K|}$.
Fixed any $i \in S$, the domain of $\overline{V} (\cdot, i)$ can be expanded to $\RR^{|K|}$ as following:
$$
\overline{\mathcal{V}} (p, i) =
\left\{
\begin{array}{ll}
\overline{V} (p, i), & p \in \pq (K); \\
-\infty, & p \in \RR^{|K|} .
\end{array}
\right.
$$
Once we prove the following equation
\be\lb{v-u1}
\inf_{z \in \RR^{|K|}} \lt \{\overline{U} (z, i) + \frac{1}{\az}\lan z, p \ran \rt\}=-\infty ,\qd \forall \, p\in \RR^{|K|}\setminus \pq (K),
\de
we can deduce that
\be\lb{v-u2}
\inf_{z \in \RR^{|K|}} \lt \{\overline{U} (z, i) + \frac{1}{\az}\lan z, p \ran \rt\}=\overline{\mathcal{V}}(p,i) ,\qd \forall \, p\in \RR^{|K|},
\de
by combining (\ref{s5-thm4-4}). Therefore, we give the proof of (\ref{v-u1}) below, which is equivalent to proof that for each $p=\{p_k,k\in K\}\in \RR^{|K|}\setminus \pq (K)$ and $M >0$, there exist some $z \in \RR^{|K|}$ such that
\be\lb{v-u3}
\overline{U} (z, i) + \frac{1}{\az}\lan z, p \ran \leq - M.
\de
If $\sum_{k\in K} p_k \neq 1$, then define a constant vector as $z = \Big\{ z(k) \equiv m_0 := \dps\frac{- \az M - c^*}{\sum_{l \in K} p_l -1}, \forall k \in K \Big\}$, and then the definition of $U$ given in (\ref{dual_eq2}) implies that
$$
\overline{U} (z, i) + \frac{1}{\az}\lan z, p \ran \leq \frac{c^*}{\az} - \frac{m_0}{\az} + \frac{m_0 \sum_{l \in K} p_l}{\az} = -M.
$$
If $\sum_{k\in K} p_k=1$ but there exists some $k_0 \in K$ satisfying $p_{k_0}<0$, then define vector $\hat{z} := \{ \hat{z} (k), \forall k \in K \}$ as
$$
\hat{z} (k) :=
\left\{
\begin{array}{ll}
z_0, & k=k_0 ; \\
0, & k \neq k_0,
\end{array}
\right.
$$
where $z_0 := (- \az M - c^*) / p_{k_0} > 0$. Again, the definition of $U$ given in (\ref{dual_eq2}) implies that
$$
\overline{U} (z, i) + \frac{1}{\az} \lan z, p \ran \leq \frac{c^*}{\az} + \sup_{q \in [0, 1]} \lt( \frac{- z_0 q}{\az} \rt) + \frac{z_0 p_{k_0}}{\az} = \frac{c^*}{\az} + \frac{z_0 p_{k_0}}{\az}= -M.
$$
Hence, it holds that (\ref{v-u3}) which guarantees that (\ref{v-u2}) holds.

According to Lemma \ref{s5-lem1} and \ref{s5-lem3}, $\overline{U} (z, i)$ is a continuous convex function on $\RR^{|K|}$. Hence, the Fenchel theorem \cite[Theorem A. 16]{SS02} implies that $\az \overline{U} (z, i) = (\az \overline{U})^{\# \#} (z, i)$, and then
\begin{align*}
\az \overline{U} (z, i) &= \sup_{p \in \RR^{|K|}} \lt\{ \lan -p, z \ran - \sup_{w \in \RR^{|K|}} \lt[ \lan w, - p \ran - (\az \overline{U}) (w, i) \rt] \rt\} \\
& = \sup_{p \in \RR^{|K|}} \lt\{ \lan -p, z \ran + \inf_{w \in \RR^{|K|}} \lt[ \lan w, p \ran + (\az \overline{U}) (w, i) \rt] \rt\} \\
& = \sup_{p \in \RR^{|K|}} \lt\{ - \lan p, z \ran + \az\overline{\mathcal{V}}(p,i) \rt\}\\
& = \sup_{p \in \pq (K)} \lt\{ \az \overline{V} (p, i) - \lan p, z \ran \rt\}.
\end{align*}
Using the compactness of $\pq(K)$ and the continuity of $\overline{V}$ again, the supremum can be replaced by the maximum, i.e., (\ref{s5-thm4-2}) holds. Similar to the case $\overline{V}$ , we expand the domain of $\underline{V} (\cdot, i)$ by
$$
\underline{\mathcal{V}} (p, i) =
\left\{
\begin{array}{ll}
\underline{V} (p, i), & p \in \pq (K); \\
-\infty, & p \in \RR^{|K|} \setminus \pq (K).
\end{array}
\right.
$$
By Lemma \ref{s3-cor6}, for each $i \in S$ we have
$$
\underline{\mathcal{V}} (\lz p_1 + (1 - \lz) p_2, i) \geq  \lz \underline{\mathcal{V}} (p_1, i) + (1- \lz) \underline{\mathcal{V}} (p_2, i), \qd \text{$\forall \lz \in [0, 1]$ and $p_1, p_2 \in \RR^{|K|}$},
$$
which means $-\underline{\mathcal{V}} (\cdot, i)$ is convex function on $\RR^{|K|}$. It can be verified that $\{ p \in \RR^{|K|} : \underline{\mathcal{V}} (p,i) \geq r \} = \{ p \in \pq(K) : \underline{V} (p,i) \geq r \}$ for each $r \in \RR$. Hence, the continuity of $V^* (\cdot, i)$ given in Lemma \ref{s3-lem4} implies that the set $\{ p \in \RR^{|K|} : \underline{ \mathcal{V}} (p,i) \geq r \}$ is closed, which means $- \underline{\mathcal{V}} (\cdot, i)$ is lower semicontinuous. According to the Fenchel theorem \cite[Theorem A. 16]{SS02}, we have
$$
- \az \underline{\mathcal{V}} (\cdot, i) = \lt(- \az \underline{\mathcal{V}} \rt)^{\# \#}(\cdot, i).
$$
Using (\ref{s5-thm4-5}) and the definition of $\underline{\mathcal{V}} $, for each $p \in \pq (K)$, we have
\begin{align*}
- \az \underline{V} (p, i) & = \sup_{z \in \RR^{|K|}} \lt\{ \lan z, p \ran - \sup_{x \in \RR^{|K|}} \lt[ \lan x, z \ran + \az\underline{\mathcal{V}} (x, i) \rt] \rt\} \\
& = \sup_{z \in \RR^{|K|}} \lt\{ \lan -z, p \ran - \sup_{x \in \pq (K)} \lt[ \lan x, -z \ran + \az \underline{V} (x, i) \rt] \rt\} \\
& = \sup_{z \in \RR^{|K|}} \lt\{ \lan -z, p \ran - \az \sup_{x \in \pq (K)} \lt[ -{\az^{-1}} \lan x, z \ran + \underline{V} (x, i) \rt] \rt\}.
\end{align*}
Combing the formula above with (\ref{s5-thm4-1}), we obtain
\begin{align*}
\underline{V} (p, i)=\inf_{z\in \RR^{|K|}} \lt\{ \underline{U}(z,i)+\frac{1}{\az} \lan p ,z\ran \rt\} \leq \min_{z\in \mathbb{B}}
\lt\{ \underline{U}(z,i)+\frac{1}{\az} \lan p ,z\ran \rt\}
\end{align*}
On the other hand, using (\ref{s5-thm4-4}) we have
$$
\underline{V} (p, i) =\overline{V}(p,i) = \min_{z\in \mathbb{B}} \lt\{ \overline{U} (z,i) + \frac{1}{\az} \lan p ,z \ran \rt\} \geq \min_{z\in \mathbb{B}} \lt\{ \underline{U} (z,i) + \frac{1}{\az} \lan p ,z\ran \rt\}.
$$
Hence, the formula (\ref{s5-thm4-3}) holds. The proof has been completed. \deprf

As the direct conclusion of Theorem \ref{s5-thm4}, we obtain the existence of value function and optimal policy of \pla 1 in the dual game $G^{\#} (z)$.

\begin{cor}\lb{s5-cor6}
Suppose that Assumption \ref{ass1} holds. For each $z \in \RR^{|K|}$ and $i \in S$, the following results hold.
\begin{description}
\item[(1)] The value function $U^* (z, i)$ of the dual game $G^{\#} (z)$ exists.
\item[(2)] There exists a optimal policy $\hat{\pz}^* \in \hat{\ppz}$ for \pla 1.
\end{description}
\end{cor}

\prf According to Theorem \ref{s3-thm3}, the value function $V^* (p, i)$ of the original game $G(p)$ exists for each $p \in \pq (K)$ and $i \in S$. For each $z \in \RR^{|K|}$ and $i \in S$, using (\ref{s5-thm4-1}) and (\ref{s5-thm4-2}), we obtain
\be\lb{s5-cor6-1}
U^* (z, i) = \overline{U} (z, i) = \underline{U} (z, i) = \max_{p \in \pq (K)} \lt\{V^* (p, i) - \frac{1}{\az} \lan p, z \ran \rt\}.
\de
The continuity of $V^* (\cdot, i)$ and the compactness of $\pq (k)$ ensure that there exists $\rz_i^* \in \pq (K)$, which depends on $i \in S$, such that
$$
U^* (z, i) = V^* (\rz_{i}^*, i) - \frac{1}{\az} \lan \rz_{i}^*, z \ran.
$$
In $G(\rho_i^*)$, using Theorem \ref{s4-thm1}, there exists a optimal policy $\pz^i =\{\pz_n^{i, (k)} , n \geq0, k\in K\}\in \ppz$ for \pla 1 such that $V (\rz_i^*, i, \pz^i, \sz) \geq V^* (\rz_i^*, i)$ for each $\sz \in \ssz$. Hence, we define $\hat{\pz}^{*} := (\rz_i^*, \pz^i)$ for each $i \in S$. Using (\ref{dual_eq1}) and (\ref{dual_eq2}), it holds that
$$
U (z, i, \hat{\pz}^*, \sz) = V (\rz^*_i, i, \pz^i, \sz) - \frac{1}{\az} \lan \rz_i^*, z \ran \geq V^* (\rz_i^*, i) - \frac{1}{\az} \lan \rz_i^*, z \ran = U^* (z, i), \qd \forall \sz \in \ssz,
$$
which means that $\hat{\pz}^{*}$ is the optimal policy for \pla 1 in $G^{\#}(z)$.
\deprf

On the whole, Theorem \ref{s5-thm4} and Corollary \ref{s5-cor6} describe the dual relations of the value function and the optimal policy for \pla 1 between $G(p)$ and $G^{\#} (z)$. Drawing on this idea, we try to use the dual game $G^{\#} (z)$ to solve the optimal policy for \pla 2 in the original $G(p)$. Previously, we introduce some natations. A set of $\RR^{|K|}$-valued bounded functions is given as
$$
\llq := \lt\{w: S \to \RR^{|K|} \ \big| \ 0 \leq w_k (i) \leq c^*, \forall i \in S, k \in K \rt\}.
$$
Since $S$ is finite, it can be verified that space $\llq$ is compact. Given any $\mz \in \ppz[0]$, $\nz \in \ssz[0]$ and $U: \RR^{|K|} \times S \to \RR$, we define a function $\ggz^{\mz, \nz} U : \pq (K) \times \llq \times \RR^{|K|} \times S \to \RR$ by
\begin{align}\lb{H}
(\ggz^{\mz, \nz} U) (p, w, z, i) &= \sum_{k} \sum_{a, b} p_k \mz^{(k)} (a | i) \nz (b | i) c (k, i, a, b) \lt( \int_0^\infty e^{-\az t} \lt( 1- D(t | i, a, b) \rt) \d t \rt) \notag \\
& \qd - \frac{1}{\az} \lan p, z \ran + \sum_{k} \sum_{a, b} p_k \mz^{(k)} (a | i) \nz (b | i) \Bigg[\sum_{j \in S} \int_0^\infty e^{-\az t} Q (\d t, j | i, a, b) \notag \\
& \qd \times \bigg( U(w(j), j) + \frac{1}{\az} \lan \llz_{\mz(\cdot | i), a} (p), w(j) \ran \bigg) \Bigg].
\end{align}
Denote by $\hat{\ppz}[0]$ the set of all functions $\phi$ satisfying
$$
\phi (k, a | i) := \rho_i (k) \mu^{(k)} (a | i), \qd \forall \rho \in \pq(K | S), \mu \in \ppz[0].
$$
It can be verified that $\hat{\ppz}[0]=\pq(K\times A|S)$. Since $K$, $S$ and $A$ are finite, the set $\hat{\ppz}[0]$ is compact and convex, see \cite{B68}. Similar to $\llz$ defined in (\ref{llz}), given any $i\in S$ and $a\in A$,  define $\chi_{(i,a)}: \hat{\ppz}[0] \to \pq(K)$ as
\be\lb{chi}
\chi_{(i,a)}(\phi)(k)=\frac{\phi(k,a|i)}{\sum_{l\in K}\phi(l,a|i)}, \qd \forall k\in K.
\de
Obviously, $\chi_{(i, a)}$ is continuous on $\hat{\ppz} [0]$. Moreover, for any $\phi_1, \phi_2 \in \hat{\ppz} [0]$ and $\lambda \in [0,1]$, it holds that
\be\lb{chi2}
\chi_{(i, a)} (\lambda \phi_1 + (1- \lambda) \phi_2) = \beta_{(i,a)} \chi_{(i,a)} (\phi_1) + (1-\beta_{(i,a)}) \chi_{(i,a)}(\phi_2),
\de
where $\beta_{(i, a)} = \dps\frac{\lambda \sum_{l \in K} \phi_1(l, a | i)}{\sum_{k\in K} [ \lambda \phi_1(k , a | i) + (1 - \lambda) \phi_2 (k, a | i)]}$.

\begin{prop}\lb{s5-lem7}
Suppose that Assumption \ref{ass1} holds. For each $z \in \RR^{|K|}$, the value function $U^* (z, i)$ of the dual game $G^{\#} (z)$ satisfies
\be\lb{s5-lem7-0}
U^* (z, i) = \min_{\nu \in \ssz[0]} \min_{w \in \llq} \max_{p \in \pq(K)} \max_{\mz \in \ppz[0]} (\ggz^{\mz, \nz} U^*) (p, w, z, i),
\de
\end{prop}

\prf According to Theorem \ref{OEthm} and Theorem \ref{s5-thm4},
\begin{align}\lb{s5-lem7-1}
& U^*(z,i) \notag \\
& = \sup_{p\in \pq(K)} \bigg\{ \sup_{\mu \in \Pi[0]} \inf_{\nu \in \ssz[0]} \sum_{k} \sum_{a, b} p_k \mz^{(k)} (a | i) \nz (b | i) c (k, i, a, b) \lt( \int_0^\infty e^{-\az t} \lt( 1- D(t | i, a, b) \rt) \d t \rt) \notag\\
& \qd + \sum_{k} \sum_{a, b} p_{k} \mu^{(k)} (a | i) \nu (b | i) \sum_{j \in \So} \int_0^{\infty} e^{-\az t} Q (\d t, j | i, a, b)V^* (\llz_{\mz (\cdot| i), a} (p), j)-\frac{1}{\az} \lan p,z \ran
\bigg\} \notag\\
& = \sup_{\phi \in \hat{\ppz}[0]}\inf_{\nz \in  \ssz[0]} H^i (\phi, \nu)
\end{align}
where
\begin{align*}
& H^i (\phi, \nu) \\
& = \sup_{\phi \in \hat{\ppz}[0]} \inf_{\nz \in \ssz[0]} \bigg\{\sum_{k} \sum_{a, b} \phi(k,a|i) \nz (b | i) c (k, i, a, b) \lt( \int_0^\infty e^{-\az t} \lt( 1- D(t | i, a, b) \rt) \d t \rt) \notag\\
& \qd + \sum_{k} \sum_{a, b} \phi(k,a|i) \nu (b | i) \sum_{j \in \So} \int_0^{\infty} e^{-\az t} Q (\d t, j | i, a, b)V^* (\chi_{(i,a)} (\phi), j) -\frac{1}{\az}\sum_{k\in K}\sum_{a\in A}\phi(k,a|i)z(k) \bigg\}.
\end{align*}
For any $\lambda \in [0,1]$ and $\phi_1, \phi_2 \in \hat{\ppz}[0]$, the concavity of $V^*$ (see Lemma \ref{s3-cor6}) and (\ref{chi2}) ensure that
\begin{align*}
\sum_{k} &\lt[ \lambda \phi_1 (k, a | i) + (1-\lambda) \phi_2 (k, a | i) \rt] V^* ( \chi_{(i,a)} (\lambda \phi_1 + (1-\lambda) \phi_2), j) \\
& \geq \sum_{k} \lt[ \lambda \phi_1(k, a | i) + (1 - \lambda) \phi_2 (k, a | i) \rt] \lt[ \beta_{(i,a)} V^* (\chi_{(i,a)} (\phi_1), j) + (1 - \beta_{i,a}) V^* (\chi_{(i,a)} (\phi_2), j) \rt] \\
& = \lambda \sum_{k} \phi_1 (k, a | i) V^* (\chi_{(i,a)} (\phi_1), j) + (1 - \lambda) \sum_{k} \phi_2 (k, a | i) V^* (\chi_{(i,a)} (\phi_2), j).
\end{align*}
This inequality implies that for any $\nu_0\in \ssz[0]$ and $r\in \RR$, the set $\hat{D} (\nu_0) := \{ \phi \in \hat{\ppz}[0] : H^i (\phi, \nu_0) \geq r\}$ is convex. According to the continuity of $\chi_{(i,a)}$ and Lemma \ref{s3-lem4}, we obtain that $\hat{D}(\nu_0)$ is closed. In a similar way, it can be showed that $\hat{E} (\phi_0) := \{\nu \in \ssz[0] : H^i (\phi_0, \nu) \leq r \}$ is convex and closed for each $\phi_0 \in \hat{\ppz[0]}$ and $r \in \RR$. Hence, the Sion minimax theorem \cite[Theorem A.7]{SS02} and the compactness of $\hat{\ppz}[0]$ and $\ssz[0]$ imply that
$$
U^* (z, i) = \min_{\nz \in \ssz[0]} \max_{\phi \in \hat{\ppz}[0]} H^i (\phi, \nu).
$$
Using the variational formula (see Theorem \ref{s5-thm4}), the value function $V^*$ in (\ref{s5-lem7-1}) can be replaced by the dual value function $U^*$, i.e.,
\begin{align*}
& U^*(z,i) \\
& = \min_{\nz \in  \ssz[0]} \max_{\phi \in \hat{\ppz}[0]} \Bigg\{ \sum_{k} \sum_{a, b} \phi(k,a|i) \nz (b | i) c (k, i, a, b) \lt( \int_0^\infty e^{-\az t} \lt( 1- D(t | i, a, b) \rt) \d t \rt) \\
& \qd + \bigg[ \sum_{k} \sum_{a, b} \phi(k,a|i) \nu (b | i) \sum_{j \in \So} \int_0^{\infty} e^{-\az t} Q (\d t, j | i, a, b) \inf_{\hat{z}\in \mathbb{B}} \Big( U^* (\hat{z}, j) + \frac{1}{\az} \lan \chi_{(i,a)} (\phi),\hat{z}\ran \Big) \bigg] \\
& \qd - \frac{1}{\az} \sum_{k} \sum_{a} \phi(k, a | i) z(k) \Bigg\} \\
& = \min_{\nz \in  \ssz[0]} \max_{\phi \in \hat{\ppz}[0]} \inf_{w \in \llq} \Bigg\{ \sum_{k} \sum_{a, b} \phi(k,a|i) \nz (b | i) c (k, i, a, b) \lt( \int_0^\infty e^{-\az t} \lt( 1- D(t | i, a, b) \rt) \d t \rt) \\
& \qd + \bigg[ \sum_{k} \sum_{a, b} \phi(k, a | i) \nu (b | i) \sum_{j \in \So} \int_0^{\infty} e^{-\az t} Q (\d t, j | i, a, b) \Big( U^* (w(j), j) + \frac{1}{\az} \lan \chi_{(i,a)} (\phi), w(j) \ran \Big) \bigg] \\
& \qd - \frac{1}{\az} \sum_{k} \sum_{a} \phi(k, a | i) z(k) \Bigg\}.
\end{align*}
According to the compactness of $\llq$ and the continuity of $U^*$ (see Lemma \ref{s5-lem1}), the infimum in the above formula can be replaced by the minimum. Since $U^*(\cdot, i)$ is convex and (see Lemma \ref{s5-lem3}), using the Sion minimax theorem \cite[Theorem A.7]{SS02}, we obtain
\begin{align*}
& U^*(z,i) \\
& = \min_{\nz \in  \ssz[0]} \min_{w \in \llq} \max_{\phi \in \hat{\ppz}[0]} \Bigg\{ \sum_{k} \sum_{a, b} \phi(k,a|i) \nz (b | i) c (k, i, a, b) \lt( \int_0^\infty e^{-\az t} \lt( 1- D (t | i, a, b) \rt) \d t \rt) \\
& \qd + \bigg[ \sum_{k} \sum_{a, b} \phi(k, a | i) \nu (b | i) \sum_{j \in \So} \int_0^{\infty} e^{-\az t} Q (\d t, j | i, a, b) \Big( U^* (w(j), j) + \frac{1}{\az} \lan \chi_{(i,a)} (\phi), w(j) \ran \Big) \bigg] \\
& \qd - \frac{1}{\az} \sum_{k} \sum_{a} \phi(k, a | i) z(k) \Bigg\}.
& = \min_{\nu \in \ssz[0]} \min_{w \in \llq} \max_{p \in \pq(K)} \max_{\mz \in \ppz[0]} (\ggz^{\mz, \nz} U^*) (p, w, z, i),
\end{align*}
which completes the proof.  \deprf

The above proposition states the optimality equation of the value function $U^*(z,i)$ in the dual game $G^{\#} (z)$, which ensures that the existence of the optimal policy for \pla 2. This is the main conclusion of this section, see Theorem \ref{s5-thm8} below.

\begin{thm}\lb{s5-thm8}
Suppose that Assumption \ref{ass1} holds. Given any $z \in \RR^{|K|}$ and $i \in S$, there exists $\sz^* \in \ssz$ such that
\be\lb{s5-thm8-0}
U (z, i, \hat{\pz}, \sz^*) \leq U^*(z, i), \qd \forall \hat{\pz} \in \hat{\ppz},
\de
i.e., $\sz^*$ is the optimal policy for \pla 2 in the dual game $G^{\#} (z)$.
\end{thm}

\prf According to Corollary \ref{s5-cor6}, the value function $U^*(z, i)$ of the dual game $G^{\#} (z)$ exists. For each $\{h_n\in \Ho_n,n\geq 0\}$ satisfying $h_0 = i$ and $h_{n}=(h_{n-1}, a_{n-1}, b_{n-1}, i_{n})$, we construct a sequence $\{ (\xz_n[h_{n-1}], \nz_n^*[h_{n}]) \in \llq \times \ssz[0], n \geq 1 \}$ by recursion. In details, $\xz_n[h_{n-1}]$ are given as
\begin{align*}
\xz_1[h_0] &={\arg\min}_{w \in \llq}\lt\{ \min_{\nz \in \ssz[0]} \max_{p \in \pq(K)} \max_{\mz \in \ppz[0]} \lt( \ggz^{\mz, \nz} U^* \rt) (p, w, z, i) \rt\}, \\
\xz_{n+1}[h_n] & = {\arg\min}_{w \in \llq} \lt\{ \min_{\nz \in \ssz[0]} \max_{p \in \pq(K)} \max_{\mz \in \ppz[0]} \lt( \ggz^{\mz, \nz} U^* \rt) (p, w, \xz_n[h_{n-1}] (i_{n}), i_{n}) \rt\}, \qd n \geq 1,
\end{align*}
and $\nz_n^*[h_{n}]$ are as
\begin{align*}
\nz^*_0[h_0] &= {\arg\min}_{\nz \in \ssz[0]}\lt\{ \min_{w \in \llq} \max_{p \in \pq(K)} \max_{\mz \in \ppz[0]} \lt( \ggz^{\mz, \nz} U^* \rt) (p, w, z, i) \rt\}, \\
\nz_n^*[h_{n}] &= {\arg\min}_{\nz \in \ssz[0]}\lt\{ \min_{w \in \llq} \max_{p \in \pq(K)} \max_{\mz \in \ppz[0]} \lt( \ggz^{\mz, \nz} U^* \rt) (p, w, \xz_n[h_{n-1}](i_n), i_n) \rt\}, \qd n \geq 1.
\end{align*}
The existence of $\nz_n^*[h_{n}]$ and $\xz_n[h_{n-1}]$ is guaranteed by the compactness and continuity, which has been discussed in detail in the proof of Proposition \ref{s5-lem7}. The definitions of $\nz_n^*[h_{n}]$ and $\xz_n[h_{n-1}]$ says that
\be\lb{s5-eq1}
U^*(z, i) = \max_{p \in \pq(K)} \max_{\mz \in \ppz[0]} \lt( \ggz^{\mz, \nz^*_0[h_0]} U^* \rt) (p, \xz_1[h_0], z, i),
\de
and for each $n \geq 1$,
\be\lb{s5-eq2}
U^* (\xz_n[h_{n-1}](i_n), i_n) = \max_{p \in \pq(K)} \max_{\mz \in \ppz[0]} \lt( \ggz^{\mz, \nz^*_n[h_n]} U^* \rt) (p, \xz_{n+1}[h_n], \xz_n[h_{n-1}](i_n), i_n).
\de
Next, we define a policy $\sz^*=\{\sz^*_n,n\geq 0\} \in \ssz$ for \pla 2, which is
$$
\sz^*_n (\cdot | h_n) := \nz_n^*[h_n] (\cdot | i_n), \qd \forall n\geq 0.
$$
The rest is to prove that $\sz^*$ is the optimal policy for \pla 2 in the dual game $G^{\#} (z)$.

To do so, we need to verify that for each $n \geq 0$ and $\hat{\pz} = (\rz, \pz) \in \hat{\ppz}$, it holds that
\begin{align}\lb{s5-thm8-1}
U^* (z, i) & \geq \frac{1}{\az} \sum_{m=0}^n \EE_{i}^{\hat{\pz}, \sz^*} \lt[ \lt(e^{-\az T_m} - e^{-\az T_{m+1}} \rt) c(\kz, X_m, A_m, B_m) \rt] - \frac{1}{\az} \lan \rz_{i}, z \ran \notag \\
& \qd + \EE^{\hat{\pz}, \sz^*}_{i} \lt[ e^{-\az T_{n+1}} \lt(U^* (\xz_{n+1}[H_n](X_{n+1}), X_{n+1}) + \frac{1}{\az} \lan _{n+1}{\rz_{i}^\pz}, \xz_{n+1}[H_n] (X_{n+1}) \ran \rt) \rt],
\end{align}
where $_n{\rz_{i}^\pz}$ is the posterior distribution on $K$ at $n$-th decision epoch. The definitions and properties of $_n{\rz_{i}^\pz}$ are given in Lemma \ref{s4-lem1}, in which the probability $p \in \pq (K)$ needs to be replaced by $\rz_i \in \pq (K)$. Hence, by Lemma \ref{s4-lem1}, we have
\be\lb{s5-thm8-2}
_n{\rz_{i}^\pz} = \llz_{\pz_{n-1} (\cdot | h_{n-1}), a_{n-1}} (_{n-1}{\rz_{i}^\pz}), \qd \forall n \geq 1.
\de
Next, we prove (\ref{s5-thm8-1}) by induction. For each $z \in \RR^{|K|}$ and $\hat{\pz} = (\rz, \pz) \in \hat{\ppz}$, since $\sz^*_0( \cdot | h_0) = \nz^*_0[h_0] ( \cdot | i)$, we have
\begin{align*}
& \EE_{i}^{\hat{\pz}, \sz^*} \lt[ \frac{1}{\az} \lt(1 -  e^{-\az T_{1}} \rt) c(\kz, X_0, A_0, B_0) \rt] - \frac{1}{\az} \lan \rz_{i}, z \ran \\
& \qd + \EE^{\hat{\pz}, \sz^*}_{i} \lt[ e^{-\az T_{1}} \lt( U^* (\xz_{1}[H_0] (X_{1}), X_{1}) + \frac{1}{\az} \lan _1{\rz_{i}^{\pz}}, \xi_1[H_0] (X_1) \ran \rt) \rt] \\
& = \sum_{k} \sum_{a, b} \rz_{i} (k) \pz_0^{(k)} (a | i) \nz_0^* [i](b | i) c (k, i, a, b) \lt( \int_0^\infty e^{-\az t} (1- D(t | i, a, b)) \d t \rt) - \frac{1}{\az} \lan \rz_{i}, z \ran \\
& \qd + \sum_{k} \sum_{a, b} \rz_{i} (k) \pz_0^{(k)} (a | i) \nz_0^*[i] (b | i) \Bigg[ \sum_{j \in S} \int_0^\infty e^{-\az t} Q (\d t, j | i, a, b) U^* (\xi_1[i] (j), j) \\
& \qd + \frac{1}{\az} \lan \llz_{\pz_0 (\cdot | i), a} (\rz_{i}), \xi_1[i] (j) \ran \Bigg] \\
& = \ggz^{\pz_0, \nz^*_0[i]} U^* ({\rz}_{i}, \xi_1[i], z, i) \\
&\leq U^* (z, i),
\end{align*}
where the last inequality is based on (\ref{s5-eq1}). That means (\ref{s5-thm8-1}) holds for $n=0$. For the case of $n \geq 1$, we need to calculate two conditional expectations. Firstly, using (\ref{dual-p1})-(\ref{dual-p3}) directly, we have
\begin{align}\lb{s5-thm8-3}
& \EE_{i}^{\hat{\pz}, \sz^*} \lt[ \frac{1}{\az} \lt( 1- e^{-\az (T_{n+2} - T_{n+1})} \rt) c (\kz, X_{n+1}, A_{n+1}, B_{n+1}) \Big| H_{n+1}, T_{m}, 0 \leq m \leq n+1 \rt] \notag \\
& = \sum_{k} \PP_{i}^{\hat{\pz}, \sz^*} \lt( \kz = k | H_{n+1}, T_{m}, 0 \leq m \leq n+1 \rt) \sum_{a, b} \pz_{n+1}^{(k)} (a | H_{n+1}) \sz_{n+1}^* (b | H_{n+1}) c(k, X_{n+1}, a, b) \notag \\
& \qd \times \int_0^{\infty} e^{-\az t} (1 - D(t | X_{n+1}, a, b)) \d t \notag \\
& = \sum_{k} {_{n+1}{\rz_{i}^\pz}} (k) \sum_{a, b} \pz_{n+1}^{(k)} (a | H_{n+1}) \nz_{n+1}^*[H_{n+1}] (b | X_{n+1}) c(k, X_{n+1}, a, b) \notag \\
& \qd \times \int_0^{\infty} e^{-\az t} (1 - D(t | X_{n+1}, a, b)) \d t.
\end{align}
Secondly, according to (\ref{s5-thm8-2}), we have
\begin{align}\lb{s5-thm8-4-1}
& \EE_{i}^{\hat{\pz}, \sz^*} \bigg[ e^{-\az (T_{n+2} - T_{n+1})} \Big( U^* (\xz_{n+2}[H_{n+1}] (X_{n+2}), X_{n+2}) \notag \\
& \qd + \frac{1}{\az} \lan _{n+2}\rz_{i}^\pz, \xz_{n+2}[H_{n+1}] (X_{n+2}) \ran \Big) \Big| H_{n+1}, T_{m}, 0 \leq m \leq n+1  \bigg] \notag \\
& = \sum_{k} \PP_{i}^{\hat{\pz}, \sz^*} \lt( \kz = k | H_{n+1}, T_{m}, 0 \leq m \leq n+1  \rt) \sum_{a, b} \pz_{n+1}^{(k)} (a | H_{n+1}) \sz_{n+1}^* (b | H_{n+1}) \notag \\
& \qd \times \bigg[ \sum_j \int_0^\infty e^{-\az t} Q(\d t, j | X_{n+1}, a, b) \bigg( U^* (\xz_{n+2}[H_{n+1}] (j), j) + \frac{1}{\az} \lan _{n+2}{\rz_{i}^\pz}, \xz_{n+2} [H_{n+1}](j) \ran \bigg) \bigg] \notag \\
& = \sum_{k} {_{n+1}}{\rz_{i}^\pz} (k) \sum_{a, b} \pz_{n+1}^{(k)} (a | H_{n+1}) \nz_{n+1}^*[H_{n+1}] (b | X_{n+1}) \bigg[ \sum_j \int_0^\infty e^{-\az t}  Q(\d t, j | X_{n+1}, a, b)  \notag \\
& \qd \times \bigg( U^* (\xz_{n+2} [H_{n+1}](j), j) + \frac{1}{\az} \lan \llz_{\pz_{n+1}(\cdot | H_{n+1}), a} (_{n+1}{\rz_{i}^\pz}), \xz_{n+2} [H_{n+1}](j) \ran \bigg) \bigg].
\end{align}
Noting that given any $(h_n,a_n,b_n)\in \Ho_n\times A \times B$, we have $\pz_{n+1}(\cdot|h_n,a_n,b_n,\cdot) \in \ppz[0]$. Using (\ref{s5-thm8-3}) and (\ref{s5-thm8-4-1}), it holds that
\begin{align}\lb{s5-thm8-4}
& \EE_{i}^{\hat{\pz}, \sz^*} \lt[ \frac{1}{\az} \lt( 1- e^{-\az (T_{n+2} - T_{n+1})} \rt) c (\kz, X_{n+1}, A_{n+1}, B_{n+1}) \Big| H_{n+1}, T_{m}, 0 \leq m \leq n+1 \rt] \notag \\
& \qd + \EE_{i}^{\hat{\pz}, \sz^*} \bigg[ e^{-\az (T_{n+2} - T_{n+1})} \Big( U^* (\xz_{n+2}[H_{n+1}] (X_{n+2}), X_{n+2}) \notag \\
& \qqd + \frac{1}{\az} \lan _{n+2}\rz_{i}^\pz, \xz_{n+2} [H_{n+1}](X_{n+2}) \ran \Big) \Big| H_{n+1}, T_{m}, 0 \leq m \leq n+1 \bigg] \notag \\
& = \ggz^{\pz_{n+1}(\cdot|H_n,A_n,B_n,\cdot), \nz_{n+1}^*[H_{n+1}]} U^* \lt( _{n+1}\rz_{i}^\pz, \xz_{n+2}[H_{n+1}], \xz_{n+1}[H_n] (X_{n+1}), X_{n+1} \rt)\notag\\
&\qd + \frac{1}{\az} \lan _{n+1}\rz_{i}^\pz, \xz_{n+1}[H_n] (X_{n+1}) \ran \notag \\
& \leq U^* (\xz_{n+1} [H_n](X_{n+1}), X_{n+1}) + \frac{1}{\az} \lan _{n+1}\rz_{i}^\pz, \xz_{n+1}[H_n] (X_{n+1}) \ran,
\end{align}
where the last inequality is based on (\ref{s5-eq2}). Hence, we calculate the conditional expectation of (\ref{s5-thm8-4}), and then
\begin{align}\lb{s5-thm8-5}
& \EE_{i}^{\hat{\pz}, \sz^*} \lt[ e^{-\az T_{n+1}} \lt(U^* (\xz_{n+1} [H_n](X_{n+1}), X_{n+1}) + \frac{1}{\az} \lan _{n+1}\rz_{i}^\pz, \xz_{n+1}[H_n] (X_{n+1}) \ran \rt) \rt] \notag \\
& \geq \EE_{i}^{\hat{\pz}, \sz^*} \lt[ \frac{1}{\az} \lt( e^{-\az T_{n+1}} - e^{-\az T_{n+2}} \rt) c (\kz, X_{n+1}, A_{n+1}, B_{n+1}) \rt] \notag \\
& \qd + \EE_{i}^{\hat{\pz}, \sz^*} \Bigg[ e^{-\az T_{n+2}} \Bigg(U^* (\xz_{n+2}[H_{n+1}](X_{n+2})), X_{n+2}) + \frac{1}{\az} \lan _{n+2}\rz_{i}^\pz, \xz_{n+2} [H_{n+1}](X_{n+2}) \ran \Bigg) \Bigg].
\end{align}
If (\ref{s5-thm8-1}) holds for some $n \geq 1$, then for the case of $n+1$, (\ref{s5-thm8-5}) implies that
\begin{align*}
U^* (z, i) & \geq \frac{1}{\az} \sum_{m=0}^{n+1} \EE_{i}^{\hat{\pz}, \sz^*} \lt[ \lt(e^{-\az T_m} - e^{-\az T_{m+1}} \rt) c(\kz, X_m, A_m, B_m) \rt] - \frac{1}{\az} \lan \rz_{i}, z \ran \notag \\
& \qd + \EE^{\hat{\pz}, \sz^*}_{i} \lt[ e^{-\az T_{n+2}} \lt(U^* (\xz_{n+2}[H_{n+1}] (X_{n+2}), X_{n+2}) + \frac{1}{\az} \lan _{n+2}{\rz_{i}^\pz}, \xz_{n+2}[H_{n+1}] (X_{n+2}) \ran \rt) \rt].
\end{align*}
Hence, (\ref{s5-thm8-1}) holds for all $n \geq 0$ by induction. Passing the limit $n \to \infty$ in (\ref{s5-thm8-1}), we obtain $U^* (z, i) \geq U (z, i, \hat{\pz}, \sz^*)$. Finally, by the arbitrariness of $\hat{\pz} \in \hat{\ppz}$ and $i \in S$, we obtain the existence of the optimal policy for \pla 2 in the dual game $G^{\#} (z)$.
\deprf

The dual game $G^{\#} (z)$ is the bridge for us to study the original game $G(p)$. There are two key points to study the existence of the optimal policy for \pla 2. One is the variational formula (Theorem \ref{s5-thm4}), the other is the existence of the optimal policy in the dual game $G^{\#} (z)$ (Theorem \ref{s5-thm8}). Back to the original game $G(p)$, we can obtain the existence of optimal control directly.

\begin{thm}\lb{s5-thm9}
Suppose that Assumption \ref{ass1} holds. Given any $p \in \pq(K)$ and $i \in S$, there exists $\sz^*\in \ssz$ such that
$$
V(p, i, \pz, \sz^*) \leq V^* (p, i), \qd \forall \pz \in \ppz, i \in S
$$
i.e., $\sz^*$ is the optimal policy for \pla 2 in $G(p)$.
\end{thm}

\prf For each $i \in S$, using Theorem \ref{s5-thm4}, there exists $z^i \in \RR^{|K|}$ such that the value function $V^*(p, i)$ of the original game $G(p)$ satisfying
$$
V^* (p, i) = U^* (z^i, i) + \frac{1}{\az} \lan p, z^i \ran.
$$
By Theorem \ref{s5-thm8}, for the vector $z^i$ given above, there exists $\sz^i=\{\sz^i_n,n\geq 0\} \in \ssz$ for \pla 2 in the dual game $G^{\#} (z^i)$ such that
$$
U (z^i, i, \hat{\pz}, \sz^i) \leq U^* (z^i, i), \qd \forall \hat{\pz} \in \hat{\ppz}.
$$
Hence, we define the policy $\sz^* =\{\sz_n^*,n\geq 0\}$ by $\sz^*_n (\cdot | h_n) = \sz^{i_0} ( \cdot | h_n)$ for each $h_0 = i_0$. Then, for arbitrary $\pz \in \ppz$, let $\hat{\pz} = (p, \pz)$, and then we have
$$
V (p, i, \pz, \sz^*) = U (z^i, i, \hat{\pz}, \sz^i) + \frac{1}{\az} \lan p, z^i \ran \leq U^* (z^i, i) + \frac{1}{\az} \lan p, z^i \ran = V^* (p, i).
$$
The arbitrariness of $\pz \in \ppz$ says that $\sz^* \in \ssz$ is the optimal policy of \pla 2 in the original game $G(p)$.
\deprf

The proof of Theorem \ref{s5-thm8} also gives the iterative algorithm of the optimal policy for \pla 2. Similar to Algorithm 1, we need to calculate $\{ (\xz_n[h_{n-1}], \nz_n^*[h_{n}]) \in \llq \times \ssz[0], n \geq 1 \}$ together. The algorithm of the optimal policy for \pla 2 is arranged in the following.

\begin{breakablealgorithm}\lb{algo-p2}
\caption{optimal policy for \pla 2} 
\hspace*{0.2in} {\bf Input:} 
The two-players zero-sum semi-Markov game with incomplete information $\{ \Ko, \So, (\Ao \times \Bo), p, Q(\cdot, \cdot|i, a, b), c(k, i, a, b) \}$; the value function $U^*$ of the dual game given by Corollary \ref{s3-Vn-prop}; for each $n \geq 0$, the history $h_n = (i_0, a_0, b_0, \ldots, i_n) \in \Ho_n$.
\begin{algorithmic}[1]
\State Compute $z^{i_0} := \arg \min_{z\in \mathbb{B}} \{U^*(z, i_0) - a^{-1} \lan p,z \ran \}$. \\
Compute
\begin{align*}
\xz_1 [h_0] & = {\arg\min}_{w \in \llq} \lt\{ \min_{\nz \in \ssz[0]} \max_{p \in \pq(K)} \max_{\mz \in \ppz[0]} \lt( \ggz^{\mz, \nz} U^* \rt) (p, w, z^{i_0}, i_0) \rt\}, \\
\nz^*_0 [h_0] &= {\arg\min}_{\nz \in \ssz[0]}\lt\{ \min_{w \in \llq} \max_{p \in \pq(K)} \max_{\mz \in \ppz[0]} \lt( \ggz^{\mz, \nz} U^* \rt) (p, w, z^{i_0}, i_0) \rt\}.
\end{align*}
\For{$n=1, 2, \ldots $} \\ 
Compute
$$
\nz_n^*[h_{n}] = {\arg\min}_{\nz \in \ssz[0]}\lt\{ \min_{w \in \llq} \max_{p \in \pq(K)} \max_{\mz \in \ppz[0]} \lt( \ggz^{\mz, \nz} U^* \rt) (p, w, \xz_n[h_{n-1}](i_n), i_n) \rt\}.
$$
Update the police at the $n$-th decision epoch $\sz^*_n (\cdot | h_n) = \nz_n^* [h_n] (\cdot | i_n)$. \\
Compute
$$
\xz_{n+1}[h_n] = {\arg\min}_{w \in \llq} \lt\{ \min_{\nz \in \ssz[0]} \max_{p \in \pq(K)} \max_{\mz \in \ppz[0]} \lt( \ggz^{\mz, \nz} U^* \rt) (p, w, \xz_n[h_{n-1}] (i_{n}), i_{n}) \rt\}.
$$
Update the sequence $\xz_{n}[h_{n-1}]$ at the $(n+1)$-th decision epoch by $\xz_{n+1}[h_n]$ and repeat step 4.
\EndFor
\end{algorithmic}
\hspace*{0.2in} {\bf Output:} 
The policy $\sz^* := \{ \sz^{*}_n, n\geq 0\}$ is the optimal policy for \pla 2.
\end{breakablealgorithm}

\bibliographystyle{alpha}

\end{document}